\documentclass[review]{elsarticle}

\usepackage{lineno,hyperref}
\modulolinenumbers[5]
\journal{XXX}

\usepackage{amsmath,amssymb,amsthm,amsfonts,afterpage,float,bm,stmaryrd,paralist,wrapfig}

\usepackage{accents}
\usepackage{cancel}
\usepackage{color}
\usepackage{fancyhdr}
\usepackage{graphics}            
\usepackage{hyperref}  
\usepackage{graphicx,epsf}
\usepackage{makeidx}
\usepackage{epsfig}
\usepackage{lscape}
\usepackage{graphicx}
\usepackage[lofdepth,lotdepth]{subfig}
\usepackage{amsmath,amssymb,amsthm,amsfonts,afterpage,float,bm,stmaryrd,paralist,wrapfig,bbm,soul}

\usepackage{cancel}
\usepackage{color}
\usepackage{fancyhdr}
\usepackage{graphics}
\usepackage{hyperref}
\usepackage{graphicx,epsf}
\usepackage{makeidx}
\usepackage{epsfig}
\usepackage{lscape}

\usepackage{listings}

\usepackage{empheq,xcolor}
\usepackage{amsmath,eqparbox,xparse,lipsum}

\pssilent

\newtheorem{theorem}{Theorem}[section]
\newtheorem{lemma}{Lemma}[section]

\newtheorem{corollary}{Corollary}[section]

\newtheorem{remark}{Remark}[section]

\numberwithin{equation}{section}
\numberwithin{figure}{section}
\numberwithin{table}{section}

\def\XXint#1#2#3{{\setbox0=\hbox{$#1{#2#3}{\int}$}
\vcenter{\hbox{$#2#3$}}\kern-.51\wd0}}

\newcommand{\cL}{\mathcal L}

\graphicspath{{Figure/}}

\begin{document}

\setlength{\pdfpageheight}{\paperheight}
\setlength{\pdfpagewidth}{\paperwidth}
\title{Asymptotically Compatible Schemes for Nonlocal
Ohta–Kawasaki Model}
\author{Wangbo Luo}
\address{Department of Mathematics, George Washington University, Washington D.C., 20052}
\author{Yanxiang Zhao\fnref{myfootnote}}
\address{Department of Mathematics, George Washington University, Washington D.C., 20052}
\fntext[myfootnote]{Corresponding author: yxzhao@email.gwu.edu}

\begin{abstract}

We study the asymptotical compatibility of the Fourier spectral method in multidimensional space for the Nonlocal Ohta-Kawasaka (NOK) model, which is proposed in our previous work \cite{Zhao_Luo1d}. By introducing the Fourier collocation discretization for the spatial variable, we show that the asymptotical compatibility holds in 2D and 3D over a periodic domain. For the temporal discretization, we adopt the second-order backward differentiation formula (BDF) method. We prove that for certain nonlocal kernels,  the proposed time discretization schemes inherit the energy dissipation law. In the numerical experiments, we verify the asymptotical compatibility, the second-order temporal convergence rate, and the energy stability of the proposed schemes. More importantly, we discover a novel square lattice pattern when certain nonlocal kernel are applied in the model. In addition, our numerical experiments confirm the existence of an upper bound for the optimal number of bubbles in 2D for some specific nonlocal kernels. Finally, we numerically explore the promotion/demotion effect induced by the nonlocal horizon $\delta$, which is consistent with the theoretical studies presented in our earlier work \cite{Zhao_Luo1d}. 
\end{abstract}


\date{\today}
\maketitle

\section{Introduction}\label{sec:Introduction}
In recent years, the Ohta-Kawasaki (OK) model, introduced in \cite{OhtaKawasaki_Macromolecules1986}, has been used for the
study of phase separation of diblock copolymers. Diblock copolymers are chain molecules composed of two distinct segment species denoted as species $A$ and $B$, respectively. These species exhibit a tendency to undergo phase separation due to their inherent chemical incompatibility, and the behavior has drawn significant attention within the field of materials science, attributed to their remarkable ability for self-assembly into nanoscale ordered structures \cite{Hamley}.

Block copolymers provide simple and easily controlled materials for exploring self-assembly phenomena. Mean field theories, accompanied by their corresponding free energy functionals, have demonstrated substantial utility in understanding and predicting pattern morphologies \cite{Bats_Fredrickson,Choksi}. Both in experimental investigations and theoretical analyses, the periodic structures from block copolymer systems have been studied over the past decades \cite{Bats_Fredrickson,Bahiana_Oono,Hasegawa_Tannaka_Yamasaki_Hashimoto,Zheng_Wang,Tang_Qiu_Zhang_Yang,Li_Jiang_Chen}. In addition to the conventional Ohta-Kawasaki theory, our focus extends to the exploration of a Nonlocal Ohta-Kawasaki (NOK) model which is characterized by a free energy functional \cite{Zhao_Luo1d}:
\begin{align}\label{functional:NOK}
E^{\text{NOK}}[u] = \int_{\Omega} \dfrac{\epsilon}{2}|(\mathcal{L}_{\delta})^{\frac{1}{2}} u|^2 + \dfrac{1}{\epsilon}W(u)\ \text{d} \mathbf{x} + \dfrac{\gamma}{2}\int_{\Omega} |(\mathcal{L}_{\delta})^{-\frac{1}{2}}(u-\omega)|^2\ \text{d}\mathbf{x},
\end{align}
with a volume constraint
\begin{align}\label{eqn:Volume}
\int_{\Omega} (u - \omega)\ \text{d}\mathbf{x} = 0.
\end{align}
Here $\Omega = \prod_{i=1}^d [-X_i, X_i] \subset \mathbb{R}^d, d =1, 2, 3$ denotes a periodic box and $0<\epsilon \ll 1$ is an interface parameter that indicates the system is in deep segregation regime. $u = u(x)$ is a phase field labeling function that indicates the density of species $A$ in the domain, and the density of species $B$ is implicitly represented by $1-u(x)$. Function $W(u) = 18(u-u^2)^2$ is a double well potential that enforces the labeling function $u(x)$ to be $0$ or $1$ in the domain. The first integral represents the oscillation-inhibiting term for the phase coarsening, promoting the formation of a larger domain, while the second integral is the oscillation-forcing term for the phase refinement, favoring multiple smaller domains. Finally, we choose the parameter $\omega\in(0, \frac{1}{2})$ to represent the volume fraction occupied by species $A$, and take $\Omega = [-\pi, \pi)^d$ for the remainder of the paper. Note that in the NOK model early proposed in  \cite{Zhao_Luo1d}, the oscillation-inhibiting term reads
\[
\int_{\Omega} \frac{\epsilon}{2}|\nabla u|^2 + \frac{1}{\epsilon} W(u)\ \text{d}\mathbf{x},
\]
which is the standard Ginzburg-Landau energy functional. In this paper, we study this term in a more general way, by replacing $|\nabla u|^2$ by a nonlocal term $|\mathcal{L}_{\delta}^{\frac{1}{2}} u|^2$. This treatment has been considered in the previous work in \cite{ Du_Yang2017, DuJuLiQiao_JCP2018}.

The non-local diffusion operator $\cL_{\delta}$ is defined as:
\begin{align}\label{eqn:nonlocal}
\cL_{\delta}u(x)=\int_{|s|\leq\delta}\rho_{\delta}(s)(u(x)-u(x+s))\text{d}s,
\end{align} 
where the kernel function $\rho_{\delta}(x)$ is a nonnegative, radial symmetric with compact support in $|s|\leq\delta$, and has a bounded second moment \cite{Du_Yang2016}. The horizon parameter $\delta>0$ is used to measure the range of nonlocal interactions introduced by $\cL_{\delta}$. The operator $\mathcal{L}_{\delta}$ is positive semi-definite. Under the given conditions on $\rho_{\delta}$, we have the weak convergence \cite{Du_nonlocalbook}  
$$
\lim_{\delta \rightarrow 0^+} \mathcal{L}_{\delta} u = \mathcal{L}_0 u := -\Delta u, \quad \forall u \in C^{\infty}(\Omega).
$$
For a comprehensive introduction to nonlocal modeling, analysis, and computation, we recommend interested readers refer to a recent monograph \cite{Du_nonlocalbook} by Du.

To study the microphase separation and pattern formation for the NOK system, we consider the $L^2$ gradient flow dynamics of the NOK model. By incorporating a penalty term into the free energy functional (\ref{functional:NOK}), we obtain an unconstrained free energy functional designed to accommodate the volume constraint (\ref{eqn:Volume}) as follows: 
\begin{align}\label{functional:pNOK}
E^{\text{pNOK}}[u] = \int_{\Omega} \dfrac{\epsilon}{2}|(\mathcal{L}_{\delta})^{\frac{1}{2}} u|^2 + \dfrac{1}{\epsilon}W(u)\ \text{d}\mathbf{x} + \dfrac{\gamma}{2}\int_{\Omega} |(\mathcal{L}_{\delta})^{-\frac{1}{2}}(u-\omega)|^2\ \text{d}\mathbf{x} + \dfrac{M}{2}\left[\int_{\Omega}u - \omega\ \text{d}\mathbf{x}\right]^2,
\end{align}
where $M\gg1$ is the penalty constant to enforce the volume constraint. Then we derive the penalized nonlocal Allen-Cahn-Ohta-Kawasaki (pNACOK) equation for the time evolution of $u(\mathbf{x},t)$ with a given initial condition $u(\mathbf{x},0)=u_0(\mathbf{x})$:
\begin{align}\label{eqn:pNOK}
    \dfrac{\partial}{\partial t} u(\mathbf{x},t) = - \epsilon\mathcal{L}_{\delta} u(\mathbf{x},t)-\dfrac{1}{\epsilon}W'(u(\mathbf{x},t))-\gamma\cL_{\delta}^{-1}(u(\mathbf{x},t)-\omega) - \mathcal{M} \int_{\Omega}( u(\mathbf{x},t) - \omega) \ \text{d}\mathbf{x},
\end{align}
with periodic boundary conditions.

For the discussion of asymptotic compatibility and energy stability in the rest of the paper, we modify the function $W(s)$ quadratically, when $|s|>M_{\text{cut}}>0$ for some $M_{\text{cut}}$.  This ensures that $W''$ has finite upper bound, which is necessary for the energy stable schemes for Ginzburg–Landau type dynamics \cite{ShenYang_DCDS2010}. We adopt the quadratic extension $\Tilde{W}(u)$ of $W(u)$ as used in \cite{ShenYang_DCDS2010} and other related citations therein. For brevity, we will continue to use $W(u)$ to represent $\Tilde{W}(u)$, and denote $ L_{W''}:= \|W''\|_{L^{\infty}}$. Here $L_{W''}$ denote the upper bound for $|W''|$. It is worth noting that for the Allen–Cahn and Cahn–Hilliard dynamics without nonlocal term, the quadratic extension for $W(s)$ guarantees that the solutions always satisfy the maximum principle \cite{CaffarelliMuler_ARMA1995}. Although the quadratic extension is unnecessary for the solutions to satisfy the maximum principle in the Allen-Cahn dynamics, it can simplify the analysis and approximation by avoiding some technical difficulties.

\subsection{Previews Work and Our Contribution}

Extensive theoretical analysis and numerical methods have been dedicated to the OK model, initially proposed by Ohka and Kawasaki in \cite{OhtaKawasaki_Macromolecules1986}. The authors in \cite{Nishiura_Ohnishi1995,Ren_Wei2000} presented a simple analogy of a binary inhibitory system derived from the OK model for diblock copolymers. Choski \cite{Choksi2012} conducted an asymptotic analysis of the OK model, establishing the existence of global minimizers. Numerical schemes have also been developed over the past years for the OK model, with a particular focus on energy-stable approaches using different gradient flow dynamics. Among these methods, the implicit midpoint spectral approximation in \cite{BenesovaMelcherSuli_SJNA2014}, the Implicit-Explicit Quadrature (IEQ) methods \cite{Yang_JCP2016,ChengYangSHen_JCP2017}, and the stabilized semi-implicit methods \cite{ShenYang_DCDS2010,XuTang_SJNA2006} have been designed for the OK model with the $H^{-1}$ gradient flow dynamics, or the Cahn–Hilliard type dynamics \cite{ShenYang_DCDS2010}. For the $L^2$ gradient flow dynamics of the OK model, both operator-splitting energy stable methods\cite{Xu_Zhao2019,Choi_Zhao2021} and maximum principle preserving methods \cite{Xu_Zhao2020} have been developed. Additionally, the Nakazawa-Ohta (NO) model, originally proposed by Nakazawa and Ohta in \cite{OhtaNakazawa_Macromolecules1993}, has garnered much attention in recent years. Ren and Wei conducted an exploration of a family of local minimizers characterized by a lamellar structure for the NO system \cite{RenWei_JNS2003}, and their subsequent work was dedicated to the pattern formation of ternary systems \cite{RenWei_ARMA2013, RenWei_ARMA2015}. Despite these efforts, a comprehensive characterization of global minimizers for the NO model remains elusive and under-researched. Gennip and Peletier focused on the one-dimensional scenario, addressing the global minimizers of the NO model in a degenerate case and under specific formulations of the long-range interaction parameters $[\gamma_{ij}]$ \cite{GennipPeletier_CVPDE2008}. On a different note, Du and Xu made a pioneering effort to systematically study the characterization of global minimizers of the NO model in non-degenerate scenarios \cite{XuDu_JNS2022}.

In recent years, extensive studies have been dedicated to explore both the mathematical and numerical aspects of nonlocal models. Notably, in \cite{DuZhou_2013,BurchLehoucq_2011,AndreuMazonRossiToledo_2010}, the authors conducted a great deal of rigorous mathematical analysis of nonlocal models. In 2012, Du et al. \cite{DuGunzburgerLehoucqZhou_SIRE2012,DuGunzburgerLehoucqZhou_M3AS2013} attempted to develop a general framework for nonlocal models. Their work involved the investigation of a class of nonlocal diffusion problems with various possible boundary conditions. Moreover, they presented numerous examples demonstrating that the nonlocal model finds various applications, ranging from continuum mechanics to graph theory. Because exact solutions for nonlocal models are rarely available, numerical methods have become crucial tools for the study of the nonlocal models. For instance, Du and Yang \cite{Du_Yang2016,Du_Yang2017} studied the Fourier spectral method for a 1D nonlocal Allen–Cahn (NAC) equation, and designed a fast and accurate implementation for the Fourier symbols of some nonlocal diffusion operators. Moreover, Du, Ju, Li, and Qiao \cite{DuJuLiQiao_JCP2018} designed stabilized linear semi-implicit schemes to study the nonlocal Cahn–Hilliard equation.  Remarkably in \cite{TianDu_SINA}, Tian and Du proposed the concept of asymptotic compatibility, to study the limiting behavior of the nonlocal model solutions $u^{\delta}_h$ as $\delta, h \rightarrow 0$.
In \cite{Du_Yang2016}, Du and Yang showed the asymptotic compatibility of the Fourier spectral collocation approximations for the NAC equation.

Recently, Zhao, Xu, and Luo \cite{Zhao_Luo1d,Xu_Zhao2020} attempt to expand the original OK model into a more general form. We introduced the NOK model (\ref{functional:NOK}) by incorporating a nonlocal operator in the NOK energy functional. Notably, our new model successfully replicated certain unconventional experimental patterns presented in \cite{XuRussellOckoChecco_SoftMatter2011,JJ_NC2020}. For the 1D case in \cite{Zhao_Luo1d}, we explore the influence of nonlocal parameters on the pattern formation for the diblock copolymer system.

The primary focus of this paper centers on the numerical analysis of the NOK model in 2D and 3D. Our contribution to this paper lies in several aspects. Firstly, we prove the asymptotic compatibility of the Fourier spectral collocation approximations applied to the pNACOK equation in 2D and 3D. The asymptotic compatibility condition is a fundamental component of the numerical analysis of the nonlocal models. Secondly, we design a second-order in-time scheme to solve the pNACOK equation, inheriting the energy dissipative law at the discrete level. Lastly, our numerical experiments reveals some unusual patterns that are completely different from those generated by the original OK model. We conduct a sequence of numerical studies on the effect of the nonlocal parameters on the bubble patterns for the NOK system. One important and novel result shows that the optimal number of bubbles in the bubble patterns may have an upper bound as the long-range interaction strength $\gamma$ goes to infinity when considering the NOK model and choosing some proper nonlocal operators in the model. This result is in contrast to the original OK model, in which the optimal number of bubbles grow to infinity as $\gamma \rightarrow \infty$. 

The rest of the paper is organized as follows. In Section \ref{sec:pre}, we first introduce the continuous (discrete) nonlocal operator $\mathcal{L}_{\delta}\ (\mathcal{L}_{\delta,h})$ along with its eigenvalues under periodic boundaries, and prove the asymptotic compatibility of Fourier spectral collocation approximation for linear nonlocal diffusion problems in 2D and 3D cases. Section \ref{sec:ac} demonstrates the asymptotic compatibility between the numerical solutions of pNACOK equations and its local counterpart in 2D and 3D. In Section \ref{sec:es}, we introduce the second-order BDF scheme for the pNACOK equation and prove its energy stability. In Section \ref{sec:ne}, a series of numerical experiments will be presented. These experiments validate the asymptotical compatibility, the rate of convergence, and the energy stability for the proposed numerical scheme. Furthermore several numerical experiments are conducted to systematically study the effect of the nonlocal parameters on the bubble patterns for the NOK system. Concluding remarks are given in Section \ref{sec:conrem}.


\section{Preliminaries and Bounds for the Fourier Symbols}\label{sec:pre}
In this section, we explore several useful bounds for the Fourier symbols (Fourier eigenvalues) for some nonlocal kernels $\mathcal{L}_{\delta}$ in 2D and 3D over periodic domain $\Omega$. The 1D results have been well established in \cite{Du_Yang2016}.
 
 \subsection{Nonlocal Diffusion Operator $\mathcal{L}_{\delta}$} \label{subsec:nonlocal}
For the nonlocal operator $\mathcal{L}_{\delta}$ in (\ref{eqn:nonlocal}), we consider radially symmetric kernel functions in $d$-dimension, $d = 1,2,3$,
 \begin{align*}
 \rho_{\delta}(|\mathbf{x}|) = \frac{1}{\delta^{d+2}}\rho\left( \frac{|\mathbf{x}|}{\delta} \right),
 \end{align*}
 where $\rho(\cdot)$ is a nonnegative nonincreasing function with a compact support in $[0,1]$.  Kernel $\rho_{\delta}$ satisfies the bounded second moment condition
\begin{align}\label{eqn:secondmoment}
\frac{1}{2}\int_{\mathbb{R}^d} \rho_{\delta}(|\mathbf{x}|) |\mathbf{x}|^2 \text{d}\mathbf{x} = d, \text{\ or\ equivalently} \ \int_0^1\rho(r)r^{d+1}\ \text{d}r = \frac{2d}{S_{d}},
\end{align}
where $S_{d}$ is the area of the unit sphere in $\mathbb{R}^d$. Then, for any $0 \neq \mathbf{k} = ( k_1,\cdots,k_d) \in \mathbb{Z}^d$, $e^{i\mathbf{k}\cdot \mathbf{x}}$ is an eigenfunction of nonlocal operator $\mathcal{L}_{\delta}$ in $\Omega = [-\pi, \pi)^d$ with the corresponding Fourier symbols \cite{Du_Yang2016,Du_Yang2017,DuTian_FCM2020}:
\begin{align}\label{eqn:eig_1d2d3d}
& \text{1D}: \  \lambda_{\delta}(\mathbf{k}) = 2\int_{0}^{\delta} \rho_{\delta}(r)\left(1-\cos(r|\mathbf{k}|)\right) \ \text{d}r,   \nonumber  \\
& \text{2D}: \  \lambda_{\delta}(\mathbf{k}) = 4 \int_0^{\frac{\pi}{2}} \int_{0}^{\delta} r \rho_{\delta}(r)\left(1-\cos(r|\mathbf{k}|\cos{\theta})\right)\ \text{d}r\text{d}\theta,     \\
& \text{3D}: \  \lambda_{\delta}(\mathbf{k}) = 4\pi\int_0^{\frac{\pi}{2}} \sin{\theta} \int_0^{\delta} r^2 \rho_{\delta}(r)\left(1-\cos(r|\mathbf{k}|\cos{\theta})\right) \ \text{d}r\text{d}\theta,   \nonumber  
\end{align}
where $|\mathbf{k}| := \|\mathbf{k}\|_2$. 

Moreover, we adopt the notations from \cite{Zhao_Xu2019,Du_Yang2016} to define several spaces. Define the space of all periodic functions in $H^s(\Omega), s\ge 0$, as $H^s_{\text{per}}(\Omega)$. Then we denote the subspace consisting of all functions $u\in H^s_{\text{per}}(\Omega)$ with zero mean as:
\begin{align}\label{eqn:space}
    \accentset{\circ}{H}_{\text{per}}^s(\Omega):=\left\{ u\in H^s_{\text{per}}(\Omega): \int_{\Omega}u(\mathbf{x})d\mathbf{x}=0    \right\}.
\end{align}
We use $\|\cdot\|_{H^s}$ to represent the standard Sobolev norm. When $s = 0$, $H^s(\Omega) = L^2(\Omega)$ and we take $\langle \cdot, \cdot \rangle$ as the $L^2$ inner product and $\|\cdot\|_{H^s} = \|\cdot\|_{L^2}$.

We define the inverse nonlocal operator $\mathcal{L}^{-1}_{\delta}:\accentset{\circ}{L}^2_{\text{per}}(\Omega)\longrightarrow\accentset{\circ}{H}^1_{\text{per}}(\Omega)$ as $\mathcal{L}_{\delta}^{-1}f=u $ if and only if $\mathcal{L}_{\delta}u=f$. In terms of Fourier series, $\mathcal{L}_{\delta}$ and $\mathcal{L}_{\delta}^{-1}$ read
\begin{align}\label{eqn:Fourier_series}
   \mathcal{L}_{\delta}u =\sum_{\mathbf{k}\in\mathbb{Z}^d\setminus{\{\mathbf{0}\}}}\lambda_{\delta}(\mathbf{k})\hat{u}_{\mathbf{k}}e^{i\mathbf{k}\mathbf{x}}, \ \ 
    \mathcal{L}_{\delta}^{-1}f=\sum_{\mathbf{k}\in\mathbb{Z}^d\setminus{\{\mathbf{0}\}}}\frac{1}{\lambda_{\delta}(\mathbf{k})}\hat{f}_{\mathbf{k}}e^{i\mathbf{k}\mathbf{x}},
\end{align}
where $\hat{u}_{\mathbf{k}}$ and $\hat{f}_{\mathbf{k}}$ are the $\mathbf{k}$-th Fourier coefficients of $u$ and $f$, respectively. 

\subsection{Discrete Nonlocal Operator $\mathcal{L}_{\delta,h}$} \label{subsection:DiscreteNonlocalOperator}

Next, we introduce some notations for the Fourier spectral collocation approximation in space. We adopt the notations from \cite{Zhao_Xu2019, DuJuLiQiao_JCP2018}. By choosing a positive even integer $N$, we take the spatial size $h =  \frac{2\pi}{N}$ and define the discrete domain $\Omega_h = \Omega\cap (\overset{d}{\underset{i = 1}{\bigotimes}} h\mathbb{Z})$, $d=1,2,3$. Two index sets related to the discrete Fourier transform are given as follows
\begin{align*}
    & S_{h} = \{\mathbf{j} = (j_1,\cdots,j_d)^T\in\mathbb{Z}^d;\ j_i = 1:N,\ i = 1:d\}, \\
    & \hat{S}_h = \{\mathbf{k} = (k_1,\cdots,k_d)^T\in\mathbb{Z}^d;\ k_l = -\frac{N}{2}+1:\frac{N}{2},\ l = 1:d\}.
\end{align*}
Denote by $\mathcal{M}_{h}$ the collection of periodic grid functions defined on $\Omega_{h}$:
\begin{align*}
    \mathcal{M}_{h} = \{f:\Omega_{h}\to \mathbb{R} | \ f_{\mathbf{j}+\mathbf{m}N} = f_{\mathbf{j}}, \ \forall\ \mathbf{j}\in S_{h}, \ \forall\ \mathbf{m} \in\mathbb{Z}^d\}.
\end{align*}
and $\accentset{\circ}{\mathcal{M}}_{h} := \{f\in\mathcal{M}_{h}|\left\langle f,1\right\rangle_{h} =0 \}$ the collection of functions in ${\mathcal{M}}_{h}$ with zero mean. For any $f,g\in \mathcal{M}_{h}$ and $\mathbf{f} = (f^1,\cdots, f^d)^T, \mathbf{g} = (g^1,\cdots, g^d)^T \in( \mathcal{M}_{h})^d$, we define the discrete $L^2$ inner product $\left\langle\cdot,\cdot\right\rangle_{h}$, discrete $L^2$-norm $\|\cdot\|_{L^2,h}$, and discrete $L^{\infty}$-norm $\|\cdot\|_{L^{\infty},h}$ as follows:
\begin{align*}
   & \langle f,g \rangle_{h} = h^d\sum_{\mathbf{j}\in S_h}f_{\mathbf{j}}g_{\mathbf{j}}, \ \ \|f\|_{L^2,h} = \sqrt{\langle f,f \rangle_{h}}, \ \ \|f\|_{L^{\infty},h} = \max_{\mathbf{j}\in S_h}|f_{\mathbf{j}}|, \\
   &  \langle \mathbf{f},\mathbf{g} \rangle_{h} =  h^d\sum_{\mathbf{j}\in S_h}\left(f^1_{\mathbf{j}}g^1_{\mathbf{j}}+\cdots +f^d_{\mathbf{j}}g^d_{\mathbf{j}}\right), \ \ ||\mathbf{f}||_{L^2,h} = \sqrt{\langle \mathbf{f},\mathbf{f} \rangle_{h}}.
\end{align*}

For a function $f\in\mathcal{M}_{h}$, we denote by $P$ the discrete Fourier transform (DFT) $\hat{f} = Pf$:
\begin{align}\label{eqn:spectral_col}
\hat{f}_{\mathbf{k}} = \frac{1}{N^d}\sum_{\mathbf{j}\in S_h} f_{\mathbf{j}} e^{-\textit{i}\mathbf{k}\cdot \mathbf{x}_{\mathbf{j}}}, \ \ \mathbf{k}\in\hat{S}_h, \ \ \mathbf{x}_{\mathbf{j}} =  -\pi \mathbf{1}_d+  h\mathbf{j},
\end{align}
where $\mathbf{1}_d$ is the $d$ dimensional all-one vector,  and the corresponding inverse discrete Fourier transform (iDFT) $f = P^{-1}\hat{f}$ is given as
\begin{align*}
{f}_{\mathbf{j}} = \sum_{\mathbf{k}\in \hat{S}_h} \hat{f}_{\mathbf{k}} e^{\textit{i}\mathbf{k} \cdot \mathbf{x}_{\mathbf{j}}}, \ \ \mathbf{j}\in{S}_h. 
\end{align*}
We have the Parseval identity holds for any function $f\in \mathcal{M}_h$, 
\begin{align}\label{eqn:plancherel_theorem}
\|f\|_{L^2,h}^2 = h^d \sum_{\mathbf{j}\in S_h} |f_{\mathbf{j}}|^2 = (2\pi)^d \sum_{\mathbf{k}\in \hat{S}_h} |\hat{f}_{\mathbf{k}}|^2.
\end{align}
By using the definition of $P$ and $P^{-1}$, we define the discrete nonlocal operator $\mathcal{L}_{\delta,h}: \mathcal{M}_{h} \rightarrow \accentset{\circ}{\mathcal{M}}_{h}$ through the use of $\mathcal{L}_{\delta,h} = P^{-1}\hat{\mathcal{L}}_{\delta,h} P$ in which 
\begin{align}
\hat{\mathcal{L}}_{\delta, h}  \hat{f}_{\mathbf{k}} := \lambda_{\delta}(\mathbf{k}) \hat{f}_{\mathbf{k}}, \quad \mathbf{k} \in \hat{S}_h.
\end{align}
In other words, for any $u \in \mathcal{M}_{h}$,
\begin{align}\label{eqn:discrete_L}
( \mathcal{L}_{\delta,h} u )_{\mathbf{j}}= (P^{-1}\hat{\mathcal{L}}_{\delta,h}P u )_{\mathbf{j}} = \sum_{\mathbf{k}\in\hat{S}_h\backslash \{\mathbf{0}\}} \hat{\mathcal{L}}_{\delta,h} \hat{u}_{\mathbf{k}} e^{i\mathbf{k}\cdot\mathbf{x}_{\mathbf{j}}} = 
\sum_{\mathbf{k}\in\hat{S}_h\backslash \{\mathbf{0}\}} \lambda_{\delta}(\mathbf{k}) \hat{u}_{\mathbf{k}} e^{i\mathbf{k}\cdot\mathbf{x}_{\mathbf{j}}}, \quad \mathbf{j} \in S_h.
\end{align}
Then we can naturally define $\mathcal{L}_{\delta,h}^{-1}: \accentset{\circ}{\mathcal{M}}_{h} \rightarrow \accentset{\circ}{\mathcal{M}}_{h}$ through the use of $\mathcal{L}_{\delta,h}^{-1} = P^{-1} \hat{\mathcal{L}}_{\delta,h}^{-1}P$ in which
\begin{align}
\hat{\mathcal{L}}_{\delta, h}^{-1}  \hat{f}_{\mathbf{k}} := \frac{1}{\lambda_{\delta}(\mathbf{k})} \hat{f}_{\mathbf{k}}, \quad \mathbf{k} \in \hat{S}_h.
\end{align}
Namely, for any $f \in \accentset{\circ}{\mathcal{M}}_{h}$,
\begin{align}\label{eqn:discrete_invL}
( \mathcal{L}_{\delta,h}^{-1} f )_{\mathbf{j}}= (P^{-1}\hat{\mathcal{L}}_{\delta,h}^{-1} P f )_{\mathbf{j}} = \sum_{\mathbf{k}\in\hat{S}_h\backslash \{\mathbf{0}\}} \hat{\mathcal{L}}_{\delta,h} \hat{f}_{\mathbf{k}} e^{i\mathbf{k}\cdot\mathbf{x}_{\mathbf{j}}} = 
\sum_{\mathbf{k}\in\hat{S}_h\backslash \{\mathbf{0}\}} \frac{1}{\lambda_{\delta}(\mathbf{k})} \hat{f}_{\mathbf{k}} e^{i\mathbf{k}\cdot\mathbf{x}_{\mathbf{j}}}, \quad \mathbf{j} \in S_h.
\end{align}
Additionally, when $\mathcal{L}_{\delta}^{-1}$ ($\mathcal{L}_{\delta,h}^{-1}$, respectively) apply to a function $f\in H^s_{\text{per}}$ ($f \in \mathcal{M}_h$, respectively) which is NOT of zero mean, we automatically take them as
\[
\mathcal{L}_{\delta}^{-1} f := \mathcal{L}_{\delta}^{-1} (f - \bar{f}), \quad \mathcal{L}_{\delta,h}^{-1} f := \mathcal{L}_{\delta,h}^{-1} (f - \bar{f}),
\]
in which $\bar{f}$ represents the mean of $f$ in both continuous and discrete sense .


 \subsection{Estimates of Some Bounds of Fourier Symbols $\lambda_{\delta}(\mathbf{k})$}
 In this subsection, we will present several estimates for the bounds of Fourier Symbols $\lambda_{\delta}(\mathbf{k})$. These estimates will be critical for proving the asymptotical compatibility and the energy stability of our numerical schemes for solving the gradient flow dynamics of the NOK model. To this end, we consider the stationary linear problems
\begin{align}\label{eqn:stationary_nonlocal_1}
\mathcal{L}_{\delta} u^{\delta}(\mathbf{x}) = f(\mathbf{x}),
\end{align}
and its local counterpart
\begin{align}\label{eqn:stationary_local_1}
\mathcal{L}_{0} u^{0}(\mathbf{x}) = f(\mathbf{x}),
\end{align}
over the periodic domain $\Omega = [-\pi, \pi)^d, d = 2, 3$. For the sake of simplicity, we assume that $u^{\delta}$ and $u^0$ are of zero mean. We use the Fourier collocation approximations defined in Subsection \ref{subsection:DiscreteNonlocalOperator} to solve the above two equations numerically, which gives
\begin{align*}
\mathcal{L}_{\delta,h} u_h^{\delta} = f_h, \quad \mathcal{L}_{0, h} u_h^{0} = f_h.
\end{align*}
Here $f_h: = f(\mathbf{x})|_{\Omega_h} \in \mathcal{M}_h$. By the definition of $\mathcal{L}_{\delta,h}^{-1}$ in (\ref{eqn:discrete_invL}), we have 
\[
u_h^{\delta} - u_h^0 = (\mathcal{L}_{\delta,h}^{-1} - \mathcal{L}_{0,h}^{-1}) f_h=P^{-1}(\hat{\mathcal{L}}_{\delta,h}^{-1} - \hat{\mathcal{L}}_{0,h}^{-1})P f_h.
\]

 \begin{lemma}\label{lemma:ACestimate}
Assuming that $u_h^{\delta}$ and $u_h^{0}$ respectively are the Fourier collocation approximate solutions of the stationary linear problems (\ref{eqn:stationary_nonlocal_1}) and (\ref{eqn:stationary_local_1}) in the periodic domain $\Omega = [-\pi,\pi)^d, d = 2, 3$, we have 
 \begin{align}\label{eqn:estimate_Lemma21}
 \|u_h^{\delta} - u_h^{0}\|_{L^2,h}\leq C\delta^2 \|f_h\|_{L^2,h},
 \end{align}\label{eqn:AC}
 where $C$ is a constant independent of both $\delta$ and $h$.
 \end{lemma}
 
 \begin{proof}

Thanks to the Parseval identity in (\ref{eqn:plancherel_theorem}), we have 
\begin{align*} 
& \| u_{h}^{\delta} - u_h^0 \|_{L^2,h}^2 = \|P^{-1}( \hat{\mathcal{L}}_{\delta,h}^{-1} - \hat{\mathcal{L}}_{0,h}^{-1})P f_h\|_{L^2,h}^2 = (2\pi)^d \sum_{\mathbf{0}\neq\mathbf{k}\in\hat{S}_h} \left| \frac{1}{\lambda_{\delta}(\mathbf{k})} - \frac{1}{\lambda_{0}(\mathbf{k})}  \right|^2 |\hat{f}_{\mathbf{k}}|^2, \\
& \|f_h\|_{L^2,h} = (2\pi)^d\sum_{\mathbf{0}\neq\mathbf{k}\in\hat{S}_h} |\hat{f}_{\mathbf{k}}|^2.
\end{align*}
 
The proof for the 1D case was provided in \cite{Du_Yang2016}. We aim to prove the asymptotic compatibility of the Fourier collocation approximations in 2D and 3D. Following the similar techniques in Lemma 1 of \cite{Du_Yang2016}, our task is simply to show the following bound
\begin{align*}
\frac{1}{\delta^2}\left|\frac{1}{\lambda_{\delta}(\mathbf{k})} - \frac{1}{\lambda_{0}(\mathbf{k})}\right| := C_{\mathbf{k}} \leq C, \ \ \forall \ 0\neq \mathbf{k} \in \mathbb{Z}^d.
\end{align*}
In the remainder of this proof, we will mainly focus on the 2D case, while the 3D case employs similar techniques.

In 2D case, the Fourier symbol $\lambda_{\delta}(\mathbf{k})$ in (\ref{eqn:eig_1d2d3d}) can be rewritten as
\begin{align*}
\lambda_{\delta}(\mathbf{k}) & = 4 \int_0^{\frac{\pi}{2}} \int_{0}^{\delta} r \rho_{\delta}(r)\left(1-\cos(r|\mathbf{k}|\cos{\theta})\right)\ \text{d}r\text{d}\theta \\
& = \frac{4}{\delta^2}\int_0^{\frac{\pi}{2}}  \int_0^{1} r \rho(r)\left(1-\cos( r \delta  |\mathbf{k}| \cos{\theta}  )\right) \ \text{d}r\text{d}\theta.
\end{align*}
Since $1 - \cos{s} \leq \frac{s^2}{2}$, we have 
\begin{align*}
\delta^2|\lambda_{\delta}(\mathbf{k})| & = 4\int_0^{\frac{\pi}{2}}  \int_0^{1} r \rho(r)\left(1-\cos( r \delta  |\mathbf{k}| \cos{\theta}  )\right) \ \text{d}r\text{d}\theta \\
& \leq 4 \int_0^{\frac{\pi}{2}} \int_0^{1} r \rho(r)\left(\frac{\delta^2 r^2 |\mathbf{k}|^2\cos^2{\theta}}{2}\right) \ \text{d}r\text{d}\theta \\
& = 4\delta^2 |\mathbf{k}|^2 \int_0^{\frac{\pi}{2}} \frac{\cos^2{\theta}}{2}\ \text{d}\theta  \int_0^{1} r^3 \rho(r) \ \text{d}r \\
& = 4\delta^2 |\mathbf{k}|^2 \frac{\pi}{8}\frac{2}{\pi} \\
& = \delta^2 |\mathbf{k}|^2 = \delta^2|\lambda_{0}(\mathbf{k})|,
\end{align*}
in which the evaluation of $\int_0^1 r^3\rho(r)\text{d}r$ is due to (\ref{eqn:secondmoment}). Using $1 - \cos{s} \geq \frac{s^2}{2} - \frac{s^4}{24}$, we also have 
\begin{align*}
\delta^2|\lambda_{\delta}(\mathbf{k})| & \geq 4 \int_0^{\frac{\pi}{2}}  \int_0^{1} r \rho(r)\left(\frac{\delta^2 r^2 |\mathbf{k}|^2\cos^2{\theta}}{2} - \frac{\delta^4 r^4 |\mathbf{k}|^4\cos^4{\theta}}{24}\right) \ \text{d}r\text{d}\theta \\
& = \delta^2 |\mathbf{k}|^2  - 4 \delta^4 |\mathbf{k}|^4 \int_0^{1} r^5 \rho(r) \text{d}r \int_0^{\frac{\pi}{2}} \frac{\cos^4{\theta}}{24}\ \text{d}\theta \\
& \ge \delta^2 |\mathbf{k}|^2  - 4\delta^4 |\mathbf{k}|^4 \int_0^{1} r^3 \rho(r) \text{d}r   \int_0^{\frac{\pi}{2}} \frac{\cos^4{\theta}}{24}\ \text{d}\theta \\ 
& = \delta^2 |\mathbf{k}|^2 - 4 \delta^4 |\mathbf{k}|^4 \frac{2}{\pi}\frac{\pi}{128} \\
& = \delta^2 |\mathbf{k}|^2 - \frac{ \delta^4 |\mathbf{k}|^4 }{16} .
\end{align*} 
Thus, for $\delta |\mathbf{k}|  \leq \pi$, it follows that
\begin{align*}
C_{\mathbf{k}}  = \frac{1}{\delta^2}\left|\frac{1}{\lambda_{\delta}(\mathbf{k})} - \frac{1}{\lambda_{0}(\mathbf{k})}\right| 
 \leq \frac{1}{\delta^2 |\mathbf{k}|^2 - \frac{\delta^4 |\mathbf{k}|^4}{16}} - \frac{1}{\delta^2 |\mathbf{k}|^2} 
 = \frac{1}{16 - \delta^2 |\mathbf{k}|^2} \leq \frac{1}{16 - \pi^2}.
\end{align*}

When $\delta |\mathbf{k}| >\pi$, the situation becomes more complicated. We need to distinguish between four cases, 
\begin{itemize}
\item Case I: $r\rho(r)$ initially decreases, reaching a minimum at $r^*\in(0,1]$;
\item Case II: $r\rho(r)$ initially increases, reaching a maximum at $r^* \in (0,1]$.
\end{itemize}

For Case I in which $r\rho(r)$ initially decreases, reaching a minimum at $r^*\in(0,1]$, we first fix $\theta\in[\frac{\pi}{3}, \cos^{-1}(\frac{1}{3}) ]$ and $\pi \le \delta|\mathbf{k}| \le \frac{\pi}{r^*}$, then we have 
\[
r^*\delta|\mathbf{k}| \cos\theta \le \frac{\pi}{2}, \quad \delta |\mathbf{k}| \cos\theta \ge \frac{\pi}{3} >1,
\]
which implies that
\begin{align}
\int_{0}^{r^*} r \rho(r) \cos{(  r \delta  |\mathbf{k}| \cos{\theta} )} \text{d}r
\le \int_{0}^{r^*} r \rho(r) \cos{(  r )} \text{d}r. 
\end{align}
Therefore when $\pi \le \delta |\mathbf{k}| \le \frac{\pi}{r^*}$, we have the bound
\begin{align}
\delta^2|\lambda_{\delta}(\mathbf{k})|
& \ge 4  \int_{\frac{\pi}{3}}^{\cos^{-1}(\frac{1}{3})}  \int_0^{r^*}  r \rho(r)\left(1-\cos(r \delta  |\mathbf{k}|\cos{\theta}  )\right) \ \text{d}r\text{d}\theta  \nonumber \\
& \ge 4  \int_{\frac{\pi}{3}}^{\cos^{-1}(\frac{1}{3})}  \int_0^{r^*}  r \rho(r)\left(1-\cos(r )\right) \ \text{d}r\text{d}\theta  = C_1 > 0, \label{eqn:bound1_CaseIII}
\end{align}
where $C_1$ is a generic constant independent of $\delta$ and $N$.

Now for $r^*\delta |\mathbf{k}| \ge \pi$, we take $\theta \in [0, \frac{\pi}{4}]$, and then $r^*\delta |\mathbf{k}| \cos\theta \ge \frac{\pi}{\sqrt{2}}$. Note that for a decreasing function $g$ over $[0,h]$, we always have
\[
\int_{0}^{h} g(r)\cos{(r)} \text{d}r \leq \int_0^{\frac{\pi}{2}} g(r) \cos(r) \text{d}r, \quad \forall h \ge  \frac{\pi}{2}.
\]
By the decrease of $r\rho(r)$ over $[0,r^*]$, and $r^*\delta|\mathbf{k}|\cos\theta \ge \frac{\pi}{\sqrt{2}} \ge \frac{\pi}{2}$, it follows that
\begin{align*}
\int_0^{r^*} r \rho(r) \cos( r\delta |\mathbf{k}| \cos\theta ) \text{d}r 
&\le \int_0^{\frac{\pi/2}{\delta |\mathbf{k}|\cos\theta}} r \rho(r) \cos( r\delta |\mathbf{k}| \cos\theta ) \text{d}r \\
&\le \int_0^{\frac{\pi/2}{\delta |\mathbf{k}|\cos\theta}} r \rho(r) \text{d}r 
\le \int_0^{\frac{r^*}{\sqrt{2}}} r \rho(r) \text{d}r.
\end{align*}
Therefore when $r^*\delta |\mathbf{k}| \ge \pi$, we have the bound
\begin{align}
\delta^2|\lambda_{\delta}(\mathbf{k})| 
& \ge 4  \int_0^{\frac{\pi}{4}}  \int_0^{r^*}  r \rho(r)\left(1-\cos(r \delta  |\mathbf{k}|\cos{\theta}  )\right) \ \text{d}r\text{d}\theta  \nonumber\\
& \geq 4 \left[\int_{0}^{\frac{\pi}{4}} \int_0^{r^*} r \rho(r) \text{d}r\text{d}\theta  - \int_0^{\frac{\pi}{4}}\int_{0}^{\frac{r^*}{\sqrt{2}}} r \rho(r)  \ \text{d}r\text{d}\theta\right], \nonumber\\
& = 4\int_0^{\frac{\pi}{4}} \int_{\frac{r^*}{\sqrt{2}}}^{r^*} r \rho(r) \ \text{d}r\text{d}\theta  = C_2 > 0, \label{eqn:bound2_CaseIII}
\end{align}
where $C_2$ is a generic constant independent of $\delta$ and $N$.

Combining the two bounds (\ref{eqn:bound1_CaseIII}) and (\ref{eqn:bound2_CaseIII}) for Case I, we get that 
\begin{align}
\text{Case I}: \quad \delta^2 |\lambda_{\delta}(\mathbf{k})| \ge \min \{C_3, C_4\} > 0,  \quad \text{for\ } \delta|\mathbf{k}| \ge \pi.
\end{align}

For Case II in which $r\rho(r)$ initially increases, reaching a maximum at $r^*\in(0,1]$, we first fix $\theta\in[\frac{\pi}{3}, \cos^{-1}(\frac{1}{3}) ]$ and $\pi \le \delta|\mathbf{k}| \le \frac{\pi}{r^*}$. Then we have 
\[
r^*\delta|\mathbf{k}| \cos\theta \le \frac{\pi}{2}, \quad \delta |\mathbf{k}| \cos\theta \ge \frac{\pi}{3} >1,
\]
which implies that
\begin{align}
\int_{0}^{r^*} r \rho(r) \cos{(  r \delta  |\mathbf{k}| \cos{\theta} )} \text{d}r
\le \int_{0}^{r^*} r \rho(r) \cos{(  r )} \text{d}r. 
\end{align}
Therefore when $\pi \le \delta |\mathbf{k}| \le \frac{\pi}{r^*}$, we have the bound,
\begin{align}
\delta^2|\lambda_{\delta}(\mathbf{k})|
& \ge 4  \int_{\frac{\pi}{3}}^{\cos^{-1}(\frac{1}{3})}  \int_0^{r^*}  r \rho(r)\left(1-\cos(r \delta  |\mathbf{k}|\cos{\theta}  )\right) \ \text{d}r\text{d}\theta  \nonumber \\
& \ge 4  \int_{\frac{\pi}{3}}^{\cos^{-1}(\frac{1}{3})}  \int_0^{r^*}  r \rho(r)\left(1-\cos(r )\right) \ \text{d}r\text{d}\theta  = C_3 > 0, \label{eqn:bound1_CaseIV}
\end{align}
where $C_3$ is a generic constant independent of $\delta$ and $N$.

Before moving to the discussion for $r^*\delta |\mathbf{k}| \ge \pi$, we need the following inequalities
\begin{align}\label{eqn:CaseII_ineq01}
\int_0^h r\rho(r) \cos(r) \text{d}r \le \int_{0}^{\frac{\pi}{2}} r\rho(r)\cos{(r)} \text{d}r, \quad \forall h\in \left [ \frac{\pi}{2}, \frac{3\pi}{2} \right ],
\end{align}
and
\begin{align}\label{eqn:CaseII_ineq02}
\int_0^{\frac{3\pi}{2}} r\rho(r) \cos(r) \text{d}r \le 0, \quad  \int_{\frac{3\pi}{2}+2n\pi}^{\frac{7\pi}{2}+2n\pi} r\rho(r)\cos{(r)} \text{d}r \leq 0, \ n = 0,1,2,\cdots.
\end{align}

Next, we choose $\theta\in[0,\frac{\pi}{4}]$ and $\pi \leq r^* \delta |\mathbf{k}| \leq \frac{3\pi}{2}$, and then $\frac{\pi}{\sqrt{2}}\leq r^* \delta |\mathbf{k}|\cos{\theta} \leq \frac{3\pi}{2}$.  Using the inequality (\ref{eqn:CaseII_ineq01}), it follows that
\begin{align}
\int_0^{r^*} r \rho(r) \cos{( r \delta  |\mathbf{k}| \cos{\theta}  )} \text{d}r & \leq \int_0^{\frac{\pi/2}{\delta |\mathbf{k}|\cos{\theta}}} r \rho(r) \cos{(r \delta  |\mathbf{k}|\cos{\theta}  )} \text{d}r \nonumber\\
 & \leq \int_0^{\frac{\pi/2}{\delta |\mathbf{k}|\cos{\theta}}} r \rho(r) \text{d}r  \leq \int_0^{\frac{r^*}{\sqrt{2}}} r \rho(r) \text{d}r. \label{estimate:n_0_CaseIV}
\end{align}

For $\theta \in [0,\frac{\pi}{4}]$ and $\pi \leq r^*\delta |\mathbf{k}| \leq\frac{7\pi}{2}+2n\pi, \ n=0, 1,2,\cdots$, we prove the existence of the upper bound of $\int_0^{r^*} r \rho(r) \cos{( r \delta  |\mathbf{k}| \cos{\theta}  )} \text{d}r$ by induction. When $n = 0$, we have $\pi \le r^*\delta |\mathbf{k}| \le \frac{7\pi}{2}$, and $\frac{\pi}{\sqrt{2}} \le  r^*\delta |\mathbf{k}| \cos\theta \le \frac{7\pi}{2}$. Note that for $r^*\delta |\mathbf{k}| \cos\theta \in [\frac{\pi}{\sqrt{2}}, \frac{3\pi}{2}]$, the upper bound is given by the estimate (\ref{estimate:n_0_CaseIV}). For $r^*\delta |\mathbf{k}| \cos\theta \in [\frac{3\pi}{2}, \frac{7\pi}{2}]$, using the first inequality in (\ref{eqn:CaseII_ineq02}), we have
\begin{align}
\int_0^{r^*} r \rho(r) \cos{( r \delta  |\mathbf{k}|\cos{\theta}  )} \text{d}r 
& \leq \int_{\frac{3\pi/2}{ \delta |\mathbf{k}|\cos{\theta}}}^{r^*} r \rho(r) \cos{(r \delta  |\mathbf{k}|\cos{\theta} )} \text{d}r \nonumber \\
 & \leq \int_{\frac{3\pi/2}{\delta |\mathbf{k}|\cos{\theta}}}^{r^*} r \rho(r)  \text{d}r \leq \int_{r^*\frac{3\pi/2}{7\pi/2}}^{r^*} r \rho(r)  \text{d}r  = \int_{\frac{3}{7}r^*}^{r^*} r \rho(r) \text{d}r. \label{estimate:n_1_CaseIV}
\end{align}
Combining the estimates (\ref{estimate:n_0_CaseIV}) and (\ref{estimate:n_1_CaseIV}), we have that for $r^*\delta |\mathbf{k}| \in [\pi, \frac{7\pi}{2}]$ or  $r^*\delta |\mathbf{k}| \cos\theta \in [\frac{\pi}{\sqrt{2}}, \frac{7\pi}{2}]$, 
\begin{align}\label{estimate:n_2_CaseIV}
\int_0^{r^*} r \rho(r) \cos{(r \delta  |\mathbf{k}|\cos{\theta}  )} \text{d}r \le \max \left\{  \int_0^{\frac{r^*}{\sqrt{2}}} r \rho(r) \text{d}r,  \int_{\frac{3}{7}r^*}^{r^*} r \rho(r) \text{d}r \right\}.
\end{align}
Assume that for $r^*\delta |\mathbf{k}| \in [\pi, \frac{7\pi}{2}+ 2(n-1)\pi ]$ or $r^*\delta |\mathbf{k}|\cos\theta \in [ \frac{\pi}{\sqrt{2}}, \frac{7\pi}{2}+ 2(n-1)\pi ]$, estimate (\ref{estimate:n_2_CaseIV}) holds. Then for the case of $n$, we consider $r^*\delta |\mathbf{k}| \in \left[ \pi, \frac{7\pi}{2}+ 2n\pi \right]$, or
$
r^*\delta |\mathbf{k}| \cos\theta  \in \left[\frac{\pi}{\sqrt{2}}, \frac{7\pi}{2}+ 2n\pi \right ].
$
For $r^*\delta |\mathbf{k}| \cos\theta  \in \left[\frac{\pi}{\sqrt{2}}, \frac{7\pi}{2}+ 2(n-1)\pi \right ]$, estimate (\ref{estimate:n_2_CaseIV}) still holds owning to the induction hypothesis. For $r^*\delta |\mathbf{k}| \cos\theta  \in \left[ \frac{3\pi}{2}+ 2n\pi, \frac{7\pi}{2}+ 2n\pi \right ]$, using the second inequality in (\ref{eqn:CaseII_ineq02}), we have that
\begin{align*}
\int_0^{r^*} r \rho(r) \cos{(r \delta  |\mathbf{k}|\cos{\theta}  )} \text{d}r 
& \leq \int_{\frac{3\pi/2+2n\pi}{\delta |\mathbf{k}|\cos{\theta}}}^{r^*} r \rho(r) \cos{(r \delta  |\mathbf{k}|\cos{\theta}  )} \text{d}r \\
 & \leq \int_{\frac{3\pi/2+2n\pi}{\delta |\mathbf{k}|\cos{\theta}}}^{r^*} r \rho(r)  \text{d}r  
 \leq \int_{r^*\frac{3\pi/2+2n\pi}{ 7\pi/2+2n\pi }}^{r^*} r \rho(r)  \text{d}r 
 \leq \int_{\frac{3}{7}r^*}^{r^*} r \rho(r) \text{d}r.
\end{align*}
Therefore, estimate (\ref{estimate:n_2_CaseIV}) holds for any $\theta\in[0, \frac{\pi}{4}]$ and $r^*\delta|\mathbf{k}|\ge\pi$. This estimate leads to 
\begin{align}
\delta^2|\lambda_{\delta}(\mathbf{k})| 
& \ge 4  \int_0^{\frac{\pi}{4}}  \int_0^{r^*}  r \rho(r)\left(1-\cos(r \delta  |\mathbf{k}|\cos{\theta}  )\right) \ \text{d}r\text{d}\theta  \nonumber \\
& \geq 4 \left[\int_{0}^{\frac{\pi}{4}} \int_0^{r^*} r \rho(r) \text{d}r\text{d}\theta  - \max \left\{ \int_0^{\frac{\pi}{4}}\int_{0}^{\frac{r^*}{\sqrt{2}}} r \rho(r)  \ \text{d}r\text{d}\theta, \int_0^{\frac{\pi}{4}}\int_{\frac{3}{7}r^*}^{r^*} r \rho(r)  \ \text{d}r\text{d}\theta   \right\}     \right], \nonumber \\
& = 4  \min \left\{   \int_0^{\frac{\pi}{4}} \int_{\frac{r^*}{\sqrt{2}}}^{r^*} r \rho(r) \ \text{d}r\text{d}\theta, \   \int_0^{\frac{\pi}{4}} \int_0^{\frac{3}{7}r^*} r \rho(r) \ \text{d}r\text{d}\theta \right\}  = C_4 > 0, \label{eqn:bound2_CaseIV}
\end{align}
where $C_4$ is a generic constant independent of $\delta$ and $N$.

Combining the two bounds (\ref{eqn:bound1_CaseIV}) and (\ref{eqn:bound2_CaseIV}) for Case II, we get that 
\begin{align}
\text{Case II}: \quad \delta^2 |\lambda_{\delta}(\mathbf{k})| \ge \min \{C_3, C_4\} > 0,  \quad \text{for\ } \delta|\mathbf{k}| \ge \pi.
\end{align}

Finally after a long discussion for Cases I and II, we obtain that
\begin{align*}
C_{\mathbf{k}} & = \frac{1}{\delta^2}\left|\frac{1}{\lambda_{\delta}(\mathbf{k})} - \frac{1}{\lambda_{0}(\mathbf{k})}\right| \leq  \frac{1}{\delta^2\lambda_{\delta}(\mathbf{k})} \leq \max\left\{\frac{1}{C_1}, \cdots, \frac{1}{C_4}\right\}, \quad \text{for}\ \delta|\mathbf{k}| \ge \pi.
\end{align*}
Therefore, by choosing $C = \max\{\frac{1}{16-\pi^2},\frac{1}{C_1}, \cdots, \frac{1}{C_4}\}$, we obtain (\ref{eqn:estimate_Lemma21}) in 2D.

For the 3D case, Fourier symbol $\lambda_{\delta}(\mathbf{k})$ in (\ref{eqn:eig_1d2d3d}) becomes
\begin{align*}
\lambda_{\delta}(\mathbf{k}) & = 4\pi \int_0^{\frac{\pi}{2}}\sin{\theta} \int_{0}^{\delta} r^2 \rho_{\delta}(r)\left(1-\cos(r|\mathbf{k}|\cos{\theta})\right)\ \text{d}r\text{d}\theta \\
& = \frac{4\pi}{\delta^2}\int_0^{\frac{\pi}{2}} \sin{\theta}  \int_0^{1} r^2 \rho(r)\left(1-\cos( r \delta  |\mathbf{k}| \cos{\theta}  )\right) \ \text{d}r\text{d}\theta.
\end{align*}
The proof in 3D employs the same tricks as that in 2D. The key differences between 3D and 2D lie in that we need to distinguish the non-increment and non-decrement of $r^2\rho(r)$ rather than $r\rho(r)$, and additionally consider the weight function $\sin\theta$, which does not introduce any new challenges. The proof is finally completed.
\end{proof}

To prove the asymptotic compatibility of the Fourier collocation approximation for the pNACOK equation later in the next section, we reformulate the result of the above lemma as a corollary, which provides several estimates related to the nonlocal operator.
\begin{corollary}\label{remark:L2bdd}
We have the following estimates in $d$-dimension, $d = 1, 2, 3$ \cite{Du_Yang2016, Zhao_Luo1d, Xu_Zhao2019}: 
\begin{align*}
& \|\mathcal{L}_{\delta,h}^{-1} -\mathcal{L}_{0,h}^{-1}\|_{L^2,h}\leq C\delta^2, \\
& \|\mathcal{L}_{0,h}^{-2}(\mathcal{L}_{\delta,h}- \mathcal{L}_{0,h})f_h\|_{L^2,h}\leq C\delta^2 \|f_h\|_{L^2,h},\\
&\|\mathcal{L}_{\delta,h}^{-1}\|_{L^2,h} \le  C\delta^2 + A, \ \delta \geq 0.
\end{align*}
where $A$ and $C$ are constants independent of $\delta$. 
\end{corollary}
Noth that the second estimate in the corollary is a direct consequence of Lemma \ref{lemma:ACestimate}. One can refer to Corollary 2 in \cite{Du_Yang2016} for the details.

 \section{Asymptotic Compatibility for the pNACOK Equations}\label{sec:ac} 
 In this section, our primary objective is to prove the asymptotic compatibility for the Fourier collocation approximate solutions of the pNACOK equation to the true solutions of the pLACOK equation. The penalized local and nonlocal ACOK equations read:
  \begin{align}
      & \text{pNACOK:}\ \  \frac{\partial u^{\delta}}{\partial t}=-\epsilon\mathcal{L}_{\delta}u^{\delta}-\frac{1}{\epsilon}W'(u^{\delta})-\gamma\mathcal{L}^{-1}_{\delta}\left(u^{\delta}-\omega\right)- M \int_{\Omega}\left[u^{\delta}-\omega\right] \text{d}\mathbf{x}, \label{eqn:pNACOK}\\
      & \text{pLACOK:}\ \ \frac{\partial u^0}{\partial t}=-\epsilon\mathcal{L}_{0}u^0-\frac{1}{\epsilon}W'(u^0)-\gamma\mathcal{L}^{-1}_{0}\left(u^0-\omega\right)- M \int_{\Omega}\left[u^0-\omega\right] \textbf{d}\mathbf{x}. \label{eqn:pLACOK}
  \end{align}
 where  $\mathcal{L}_{0} = -\Delta$ and $\mathcal{L}_0^{-1} = (-\Delta)^{-1}$ represent the local operator and its inverse, and $u^{\delta}$ and $u^0$ denote the true solutions of the pNACOK and pLACOK equations, respectively. We denote by $u^{\delta}_h$ and $u^0_h$ the Fourier collocation approximate solutions of the pNACOK and pLACOK equations, respectively, which satisfy  \begin{align}
     & \frac{\partial u^{\delta}_h}{\partial t}=-\epsilon\mathcal{L}_{\delta,h}u^{\delta}_h-\frac{1}{\epsilon} W'(u^{\delta}_h) -\gamma \mathcal{L}^{-1}_{\delta,h}(u^{\delta}_h-\omega) - M \left\langle u^{\delta}_h-\omega,1\right\rangle_h,      \label{eqn:pNACOK_spatial}\\
     & \frac{\partial u^{0}_h}{\partial t}=-\epsilon\mathcal{L}_{0,h}u^{0}_h-\frac{1}{\epsilon} W'(u^{0}_h) -\gamma \mathcal{L}^{-1}_{0,h}(u^{0}_h-\omega) - M \left\langle u^{0}_h - \omega,1\right\rangle_h. \label{eqn:pLACOK_spatial}
 \end{align}
  
Building upon the framework outlined in \cite{TianDu_SINA}, we say that a numerical method for some nonlocal model is asymptotically compatible if
 \begin{align*}
     \|u^{\delta}_h-u^0\|_{L^2,h}\to0 \quad \text{as} \ \ \delta\to0, \ N \to \infty.
 \end{align*}
 To study the asymptotic compatibility for the pNACOK model, we apply the triangle inequality
 \begin{align*}
     \|u^{\delta}_h-u^0\|_{L^2,h}\leq\|u^{\delta}_h-u^0_h\|_{L^2,h}+\|u^{0}_h-u^0\|_{L^2,h},
 \end{align*}
 and aim to show that respectively $\|u^{\delta}_h-u^0_h\|_{L^2,h}$, $\|u^{0}_h-u^0\|_{L^2,h}\to0$ as $\delta\to0$, $N \to \infty$. 
 
 The convergence $\|u^{0}_h-u^0\|_{L^2,h}\to0$ has been established in \cite{Zhao_Xu2019}. In this work, the authors apply the Fourier collocation approximation in space to study the pLACOK equation and have the error estimate:
 \begin{align}\label{existfourierresult}
    \|u^{0}_h-u^0\|_{L^2}\leq CN^{-m},
\end{align}
where $u^0 \in H^{m}(\Omega)$ and $C$ is constant independent of $N$. Therefore, it suffices to verify that
 \begin{align}\label{eqn:ACerror}
     \|u^{\delta}_h-u^0_h\|_{L^2,h}\to0, \quad \text{as} \ \ \delta\to0, \ N \to \infty.
 \end{align}
 
To this end, we first need a lemma for the $L^2$ bound of $u^{\delta}_h$.
\begin{lemma}[$L^2$ bound for $u^{\delta}_h$]\label{lemma:u_delta}
Assuming that $u_h^{\delta}$ is the Fourier collocation approximate solution of the pNACOK equation (\ref{eqn:pNACOK_spatial}), we have 
 \begin{align}\label{eqn:L^2estimate}
 \|u_h^{\delta} - \omega\|_{L^2,h}\leq C,
 \end{align}
 where $C$ is a constant independent of both $\delta$ and $h$.
\end{lemma}
\begin{proof}
The equation (\ref{eqn:pNACOK_spatial}) can be reformulated as 
\begin{align}\label{eqn:pNACOK_newform}
\frac{\partial (u^{\delta}_h - \omega)}{\partial t} = -\epsilon\mathcal{L}_{\delta,h}(u^{\delta}_h-\omega) -\frac{1}{\epsilon} W'(u^{\delta}_h) -\gamma \mathcal{L}^{-1}_{\delta,h} (u^{\delta}_h-\omega) - M \left\langle u^{\delta}_h - \omega,1\right\rangle_h
\end{align} 
Defining the function $F(t) = \|u^{\delta}_h(t,\cdot)-\omega\|_{L^2,h}$, we have 
\begin{align*}
\frac{1}{2}\frac{d}{dt}F^2= F\frac{dF}{dt} = \left\langle \frac{\partial (u^{\delta}_h - \omega)}{\partial t},  u^{\delta}_h - \omega \right\rangle_h.
\end{align*}
Using equation (\ref{eqn:pNACOK_newform}), it follows that
\begin{align*}
F\frac{dF}{dt} & = -\epsilon\left\|\mathcal{L}_{\delta,h}^{\frac{1}{2}}(u^{\delta}_h-\omega)\right\|_{L^2,h}^2 - \frac{1}{\epsilon}\left\langle W'(u^{\delta}_h),u^{\delta}_h-\omega\right\rangle_h - \gamma\left\| \mathcal{L}^{-\frac{1}{2}}_{\delta,h}(u^{\delta}_h - \omega)\right\|_{L^2,h}^2 - \mathcal{M}\left\langle u^{\delta}_h-\omega,1 \right\rangle_h^2 \\
& \leq - \frac{1}{\epsilon}\left\langle W'(u^{\delta}_h)-W'(\omega)+W'(\omega),u^{\delta}_h-\omega\right\rangle_h\\
& \leq  \frac{L_{W''}}{\epsilon}\|u^{\delta}_h-\omega\|_{L^2,h}^2 + \frac{W'(\omega)}{\epsilon}\|u^{\delta}_h-\omega\|_{L^2,h} \\
& = \frac{L_{W''}}{\epsilon}F^2 + \frac{W'(\omega)}{\epsilon}F,
\end{align*}
which yields 
\begin{align*}
\frac{dF}{dt} \leq \frac{L_{W''}}{\epsilon}F + \frac{W'(\omega)}{\epsilon}.
\end{align*}
Finally, using the Gronwall’s inequality, we conclude that 
\begin{align*}
F =  \|u_h^{\delta} - \omega\|_{L^2,h}  \leq C, 
\end{align*}
where $C$ is a generic constant independent of both $\delta$ and $N$.
\end{proof}

Now we are ready to present the result regarding the convergence (\ref{eqn:ACerror}).

  \begin{lemma}\label{acdelta}
  Assume that $u_h^{\delta}$ and $u_h^0$ are the Fourier collocation approximate solutions of the pNACOK equation (\ref{eqn:pNACOK}) and the pLACOK equation (\ref{eqn:pLACOK}), then we have
  \begin{align*}
      \|u^{\delta}_h-u^0_h\|_{L^2,h}\leq C\delta^2,
  \end{align*}
  where $C$ is a generic constant independent of $\delta$ and $N$.
  \end{lemma}

\begin{proof}

 Following the similar method from \cite{Du_Yang2016}, define the error function
 \begin{align*}
     E(t)=\|u^{\delta}_h(t,\cdot)-u^0_h(t,\cdot)\|_{L^2,h},
 \end{align*}
 then
 \begin{align*}
     E^2(t)=\|u^{\delta}_h-u^0_h\|_{L^2,h}^2=\left\langle u^{\delta}_h-u^0_h,u^{\delta}_h-u^0_h\right\rangle_{h}.
 \end{align*}
Taking the derivative with respect to $t$, it becomes: 
 \begin{align*}
     \frac{1}{2}\frac{d}{dt}E^2 = E\frac{dE}{dt} = \left\langle \frac{\partial}{\partial t}u^{\delta}_h-\frac{\partial}{\partial t}u^0_h,u^{\delta}_h-u^0_h\right\rangle_h.
 \end{align*}
 Using the equations (\ref{eqn:pNACOK_spatial}) and (\ref{eqn:pLACOK_spatial}), we have
  \begin{align*}
      E\frac{dE}{dt}=
      & \underbrace{-\epsilon\left\langle\mathcal{L}_{\delta,h}u^{\delta}_h-\mathcal{L}_{0,h}u^{0}_h,u^{\delta}_h-u^{0}_h\right\rangle_h}_{\text{I}}\underbrace{-\frac{1}{\epsilon}\left\langle \left(W'(u^{\delta}_h)-W'(u^{0}_h)\right),u^{\delta}_h-u^{0}_h\right\rangle_h}_{\text{II}} \\
      & \underbrace{-\gamma\left\langle \left(\mathcal{L}^{-1}_{\delta,h}[u^{\delta}_h - \omega]-\mathcal{L}^{-1}_{0,h}[u^{0}_h - \omega]\right), u^{\delta}_h-u^{0}_h\right\rangle_h}_{\text{III}}  \underbrace{-\mathcal{M}\left\langle u^{\delta}_h-u^{0}_h,1 \right\rangle_h\left\langle 1,u^{\delta}_h-u^{0}_h\right\rangle_h}_{\text{IV}}.
  \end{align*}
  
 For term I, since the discrete nonlocal operator $\mathcal{L}_{\delta,h}$ is positive semi-definite,
  \begin{align*}
      \text{I} 
      & = -\epsilon\left\langle\mathcal{L}_{\delta,h}u^{\delta}_h-\mathcal{L}_{\delta,h}u^{0}_h +  \mathcal{L}_{\delta,h}u^{0}_h -\mathcal{L}_{0,h}u^{0}_h,u^{\delta}_h-u^{0}_h\right\rangle_h \\
      & = \epsilon \left\langle - \mathcal{L}_{\delta,h}[u^{\delta}_h - u^{0}_h],u^{\delta}_h-u^{0}_h\right\rangle_h - \epsilon \left\langle (\mathcal{L}_{\delta,h} - \mathcal{L}_{0,h})u^{0}_h, u^{\delta}_h-u^{0}_h \right\rangle_h \\
     &  \leq \epsilon \|(\mathcal{L}_{\delta,h} - \mathcal{L}_{0,h})u^{0}_h\|_{L^2,h}\|u^{\delta}_h-u^{0}_h\|_{L^2,h} \\
     & = \epsilon \|(\mathcal{L}_{\delta,h} - \mathcal{L}_{0,h})u^{0}_h\|_{L^2,h}E,
  \end{align*}
Note that $\|(\mathcal{L}_{\delta,h} - \mathcal{L}_{0,h})u^{0}_h\|_{L^2,h} = \|\mathcal{L}_{0,h}^{-2}(\mathcal{L}_{\delta,h} - \mathcal{L}_{0,h})\mathcal{L}_{0,h}^2u^{0}_h\|_{L^2,h}$,  we can use Corollary \ref{remark:L2bdd} to get
\begin{align*}
\text{I} \le  \epsilon \|\mathcal{L}_{0,h}^{-2}(\mathcal{L}_{\delta,h} - \mathcal{L}_{0,h})\mathcal{L}_{0,h}^2u^{0}_h\|_{L^2,h}  \le \epsilon C \delta^2 \|\mathcal{L}^{2}_{0,h}u^{0}_h\|_{L^2,h} E,
\end{align*}
where $\|\mathcal{L}^{2}_{0,h}u^{0}_h\|_{L^2,h}$ is bounded independent of $\delta$ and $N$ by Theorem 4 in \cite{Du_Yang2016} .

For term II, consider the Taylor Expansion
  \begin{align*}
      W'(u^{\delta}_h)=W'(u^0_h)+(u^{\delta}_h-u^0_h)W''(\xi),
  \end{align*}
  where $\xi$ is between $u^{\delta}_h$ and $u^0_h$,
  then
  \begin{align*}
      \text{II}
       =-\frac{1}{\epsilon}\left\langle (u^{\delta}_h-u^0_h)W''(\xi)),u^{\delta}_h-u^{0}_h\right\rangle_h 
       \le \frac{L_{W''}}{\epsilon}E^2.
  \end{align*}

  For term III, we have that
  \begin{align*}
      & \left\langle \left(\mathcal{L}^{-1}_{\delta,h}[u^{\delta}_h - \omega ]-\mathcal{L}^{-1}_{0,h}[u^{0}_h -\omega ]\right), u^{\delta}_h-u^{0}_h\right\rangle_h \\
      = & \underbrace{\left\langle (\mathcal{L}^{-1}_{\delta,h}-\mathcal{L}^{-1}_{0,h})[u^{\delta}_h-\omega],u^{\delta}_h-u^{0}_h\right\rangle_h }_{\text{a}}+\underbrace{\left\langle \mathcal{L}^{-1}_{0,h}[u^{\delta}_h-u^{0}_h],u^{\delta}_h-u^{0}_h\right\rangle_h}_{\text{b}} .
  \end{align*}
  From Corollary \ref{remark:L2bdd}, Lemma \ref{lemma:u_delta}, we have 
  \begin{align*}
      |\text{III}(\text{a})|
      & \leq \| (\mathcal{L}^{-1}_{\delta,h}-\mathcal{L}^{-1}_{0,h})[u^{\delta}_h-\omega]\|_{L^{2},h}\|u^{\delta}_h-u^{0}_h\|_{L^2,h} \leq C\delta^2E, \\
      |\text{III}(\text{b})|
      & \leq \|\mathcal{L}^{-1}_{0,h}[u^{\delta}_h-u^{0}_h]\|_{L^2,h}\|u^{\delta}_h-u^{0}_h\|_{L^2,h}  = \|\mathcal{L}^{-1}_{0,h}\|_{L^2,h} E^2.
  \end{align*}
Thus, the estimate for term III follows: 
  \begin{align*}
      \text{III} \leq \gamma \|\mathcal{L}^{-1}_{0,h}\|_{L^{2},h} E^2 + \gamma C\delta^2E .
  \end{align*}

For term IV, we have
 \begin{align*}
 \text{IV} 
 \leq  \mathcal{M}\left\langle u^{\delta}_h-u^{0}_h,1 \right\rangle_h^2 \leq \mathcal{M} |\Omega| \|u^{\delta}_h-u^{0}_h\|_{L^2,h}^2  = \mathcal{M} |\Omega| E^2, \quad |\Omega| = (2\pi)^d.
 \end{align*}

Combining I - IV, it follows  
  \begin{align*}
      E\frac{dE}{dt}\leq &  (C_1+C_2)\delta^2 E + C_3 E^2.
  \end{align*}
  Here
  \[
  C_1 = \epsilon C\|\mathcal{L}^{2}_{0,h}u^{0}_N\|_{L^2,h}, \quad C_2 = \gamma C, \quad C_3 = \frac{L_{W''}}{\epsilon} + \gamma\|\mathcal{L}^{-1}_{0,h}\|_{L^{2},h}+\mathcal{M} |\Omega|
  \]
are all independent of $\delta$ and $N$. Then the inequality becomes
  \begin{align}\label{ieqn:energy_ac}
      \frac{dE}{dt}\leq C_3E+(C_1+C_2)\delta^2.
  \end{align}
  Applying Grownwall's Inequality to (\ref{ieqn:energy_ac}), that leads to
  \begin{align*}
      E\leq C_4\delta^2,
  \end{align*}
  where $C_4$ is a constant independent of $\delta$ and $N$.
  
  \end{proof}
  
  Then, the existing result (\ref{existfourierresult}) and Lemma \ref{acdelta} imply the asymptotic compatibility.
  \begin{theorem}\label{thm:ac}
   Assume that $u_h^{\delta}$ is the Fourier collocation approximate solution to the pNACOK equation (\ref{eqn:pNACOK_spatial}), $u^0 \in H^m(\Omega)$ is the exact solution to the pLACOK equation (\ref{eqn:pLACOK}), then we have 
  \begin{align*}
      \|u^{\delta}_h-u^0\|_{L^2}\leq C (\delta^2+N^{-m} ),
  \end{align*}
  where $C$ is a generic constant independent of $\delta$ and $N$.
  \end{theorem}

\section{Second Order Time-discrete Energy Stable Schemes}\label{sec:es}

In this section, we consider the second order Backward Differentiation Formula (DBF) discretization in time for the pNACOK equation (\ref{eqn:pNOK}). With the space being discretized by Fourier collocation approximation, we now have a fully-discrete scheme. Subsequently, we will explore the energy stability at the fully-discrete level. The energy stability analysis has been extensively studied in the past decade, so we will only verify it briefly here.

Given a time interval $[0,T]$ and an integer $N>0$, we take the uniform time step size $\tau = \frac{T}{N}$ and $t_n = n\tau$ for $n=0,1,\cdots,N$. we denote by $ U^n_{\delta} \in \mathcal{M}_h:  (U^n_{\delta})_{\mathbf{j}} \approx u^{\delta}(\mathbf{x}_{\mathbf{j}};t_n)$ the approximate solution at $\mathbf{x}_{\mathbf{j}}\in \Omega_h$ and time $t_n$. Given initial conditions $U^{-1}_{\delta},  U^{0}_{\delta}$, we aim to seek $U_{\delta}^{n+1}  \in \mathcal{M}_h$ such that
\begin{align}\label{eqn:BDF_NOK_fully}
    \frac{3U^{n+1}_{\delta}-4U^{n}_{\delta}+U^{n-1}_{\delta}}{2\tau} = & -\epsilon\mathcal{L}_{\delta,h} U^{n+1}_{\delta} - \frac{1}{\epsilon}\left[2W'(U^n_{\delta})-W'(U^{n-1}_{\delta})\right] \nonumber \\
    & -\gamma B_h\mathcal{L}^{-1}_{\delta,h}\left(U_{\delta}^{n+1}-2U_{\delta}^{n}+U_{\delta}^{n-1}\right)-A_h\left(U_{\delta}^{n+1}-2U_{\delta}^{n}+U_{\delta}^{n-1}\right) \nonumber\\
    & -\gamma\mathcal{L}^{-1}_{\delta,h}\left[2U_{\delta}^n-U_{\delta}^{n-1}-\omega\right] -M\left[\langle2U_{\delta}^n-U_{\delta}^{n-1},1\rangle_h - \omega|\Omega|\right],
 \end{align}
where $A_h$, $B_h$ are stabilization constants, the stabilizer $A_h\left(U_{\delta}^{n+1}-2U_{\delta}^{n}+U_{\delta}^{n-1}\right)$ controls the growth of $W'$, and the stabilizer $\gamma B_h\mathcal{L}^{-1}_{\delta,h}\left(U_{\delta}^{n+1}-2U_{\delta}^{n}+U_{\delta}^{n-1}\right)$ dominates the behavior of $\mathcal{L}^{-1}_{\delta,h}$.

Now we briefly verify the energy stability for the proposed scheme. Taking $L^2$ inner product with respect to $U_{\delta}^{n+1}-U_{\delta}^n$ on the two sides of (\ref{eqn:BDF_NOK_fully}) yields
\begin{align}\label{eqn:LHS1}
   \text{ LHS}=&\ \frac{1}{2\tau}\left\langle2(U_{\delta}^{n+1}-U_{\delta}^n), U_{\delta}^{n+1}-U_{\delta}^n\right\rangle_h+\frac{1}{2\tau}\left\langle U_{\delta}^{n+1}-2U_{\delta}^{n}+U_{\delta}^{n-1},U_{\delta}^{n+1}-U_{\delta}^n\right\rangle_h, \\
    \text{RHS}=&\ \underbrace{-\epsilon\left\langle\mathcal{L}_{\delta,h}U_{\delta}^{n+1},U_{\delta}^{n+1}-U_{\delta}^{n}\right\rangle_h}_{\text{I}}\underbrace{-\frac{1}{\epsilon}\left\langle 2W'(U_{\delta}^{n})-W'(U_{\delta}^{n-1}),U_{\delta}^{n+1}-U_{\delta}^{n}\right\rangle_h}_{\text{II}} \nonumber\\
    & \underbrace{-A_h\left\langle U_{\delta}^{n+1}-2U_{\delta}^{n}+U_{\delta}^{n-1},U_{\delta}^{n+1}-U_{\delta}^{n}\right\rangle_h}_{\text{III}}\underbrace{-\gamma B_h\left\langle \mathcal{L}^{-1}_{\delta,h}(U_{\delta}^{n+1}-2U_{\delta}^{n}+U_{\delta}^{n-1}),U_{\delta}^{n+1}-U_{\delta}^{n}\right\rangle_h}_{\text{IV}} \nonumber \\
    & \underbrace{-\gamma\left\langle\mathcal{L}^{-1}_{\delta,h}[2U_{\delta}^{n}-U_{\delta}^{n-1}-\omega],U_{\delta}^{n+1}-U_{\delta}^{n}\right\rangle_h}_{\text{V}} \underbrace{-M\left(\left\langle 2U_{\delta}^n-U_{\delta}^{n-1},1\right\rangle_h - \omega|\Omega|\right)\left\langle 1,U_{\delta}^{n+1}-U_{\delta}^{n}\right\rangle_h}_{\text{VI}}. \nonumber
\end{align}
Then, using the identities $a\cdot(a-b)=\frac{1}{2}a^2-\frac{1}{2}b^2+\frac{1}{2}(a-b)^2$ and $b\cdot(a-b)=\frac{1}{2}a^2-\frac{1}{2}b^2-\frac{1}{2}(a-b)^2$, with $a=U_{\delta}^{n+1}-U_{\delta}^{n}$, $b=U_{\delta}^{n}-U_{\delta}^{n-1}$, (\ref{eqn:LHS1}) can be reformulated as 
\begin{align*}
    \text{LHS}=\frac{1}{\tau}\|U_{\delta}^{n+1}-U_{\delta}^{n}\|_{L^2,h}^2+\frac{1}{4\tau}\left(\|U_{\delta}^{n+1}-U_{\delta}^{n}\|_{L^2,h}^2-\|U_{\delta}^{n}-U_{\delta}^{n-1}\|_{L^2}^2+\|U_{\delta}^{n+1}-2U_{\delta}^{n}+U_{\delta}^{n-1}\|_{L^2,h}^2\right).
\end{align*}
Term I becomes: 
\begin{align*}
    \text{I}=-\frac{\epsilon}{2}\left\langle\mathcal{L}_{\delta,h}U_{\delta}^{n+1},U_{\delta}^{n+1}\right\rangle_h+\frac{\epsilon}{2}\left\langle\mathcal{L}_{\delta,h}U_{\delta}^{n},U_{\delta}^{n}\right\rangle_h-\frac{\epsilon}{2}\left\langle\mathcal{L}_{\delta,h}(U_{\delta}^{n+1}-U_{\delta}^{n}),U_{\delta}^{n+1}-U_{\delta}^{n}\right\rangle_h.
\end{align*}
For term II, by Taylor Expansion, we have
\begin{align*}
    & W'(U_{\delta}^{n})(U_{\delta}^{n+1}-U_{\delta}^{n})=W(U_{\delta}^{n+1})-W(U_{\delta}^{n})-\frac{W''(\xi_h^n)}{2}(U_{\delta}^{n+1}-U_{\delta}^{n})^2, \\
   &  W'(U_{\delta}^{n})-W'(U_{\delta}^{n-1})=W''(\mu_h^n)(U_{\delta}^{n}-U_{\delta}^{n-1}),
\end{align*}
where $\xi_h^n$ is between $U_{\delta}^{n+1}$ and $U_{\delta}^n$, $\mu_h^n$ is between $U_{\delta}^{n-1}$ and $U_{\delta}^n$. Then term II becomes
\begin{align*}
    \text{II} & =-\frac{1}{\epsilon}\left\langle 1,W'(U_{\delta}^n)(U_{\delta}^{n+1}-U_{\delta}^n)\right\rangle_h-\frac{1}{\epsilon}\left\langle W'(U_{\delta}^n)-W'(U_{\delta}^{n-1}),U_{\delta}^{n+1}-U_{\delta}^n\right\rangle_h\\
    & =-\frac{1}{\epsilon}\left\langle 1,W(U_{\delta}^{n+1})\right\rangle_h+\frac{1}{\epsilon}\left\langle 1,W(U_{\delta}^{n})\right\rangle_h+\frac{W''(\xi_h^n)}{2\epsilon}\|U_{\delta}^{n+1}-U_{\delta}^n\|_{L^2,h}^2-\frac{W''(\mu_h^n)}{\epsilon}\left\langle U_{\delta}^n-U_{\delta}^{n-1},U_{\delta}^{n+1}-U_{\delta}^n\right\rangle_h.
\end{align*}
Moreover, terms III and IV becomes
\begin{align*}
    \text{III}= &-\frac{A_h}{2}\|U_{\delta}^{n+1}-U_{\delta}^n\|_{L^2,h}^2+\frac{A_h}{2}\|U_{\delta}^{n}-U_{\delta}^{n-1}\|_{L^2,h}^2-\frac{A_h}{2}\|U_{\delta}^{n+1}-2U_{\delta}^n+U_{\delta}^{n-1}\|_{L^2,h}^2, \\
    \text{IV}=& -\frac{\gamma B_h}{2}\left\langle\mathcal{L}^{-1}_{\delta,h}(U_{\delta}^{n+1}-U_{\delta}^{n}),U_{\delta}^{n+1}-U_{\delta}^n\right\rangle_h+\frac{\gamma B_h}{2}\left\langle\mathcal{L}^{-1}_{\delta,h}(U_{\delta}^{n}-U_{\delta}^{n-1}),U_{\delta}^{n}-U_{\delta}^{n-1}\right\rangle_h \\
    & -\frac{\gamma B_h}{2}\left\langle\mathcal{L}^{-1}_{\delta,h}(U_{\delta}^{n+1}-2U_{\delta}^{n}+U_{\delta}^{n-1}),U_{\delta}^{n+1}-2U_{\delta}^{n}+U_{\delta}^{n-1}\right\rangle_h,
\end{align*}
and terms V and VI change into
\begin{align*}
    \text{V} = & -\gamma\left\langle\mathcal{L}^{-1}_{\delta,h}[U_{\delta}^n-\omega],U_{\delta}^{n+1}-U_{\delta}^n\right\rangle_h-\gamma\left\langle\mathcal{L}^{-1}_{\delta,h}[U_{\delta}^n-U_{\delta}^{n-1}],U_{\delta}^{n+1}-U_{\delta}^n\right\rangle_h \\
     = & -\frac{\gamma}{2}\left\langle\mathcal{L}^{-1}_{\delta,h}[U_{\delta}^{n+1}-\omega],U_{\delta}^{n+1}-\omega\right\rangle_h+\frac{\gamma}{2}\left\langle\mathcal{L}^{-1}_{\delta,h}[U_{\delta}^{n}-\omega],U_{\delta}^{n}-\omega\right\rangle_h \\
     & +\frac{\gamma}{2}\left\langle\mathcal{L}^{-1}_{\delta,h}[U_{\delta}^{n+1}-U_{\delta}^n],U_{\delta}^{n+1}-U_{\delta}^n\right\rangle_h - \gamma\left\langle\mathcal{L}^{-1}_{\delta,h}[U_{\delta}^n-U_{\delta}^{n-1}],U_{\delta}^{n+1}-U_{\delta}^n\right\rangle_h,\\
    \text{VI} = & -M \left\langle U_{\delta}^n - \omega,1\right\rangle_h   \left\langle 1,U_{\delta}^{n+1}-U_{\delta}^{n}\right\rangle_h -M\left\langle U_{\delta}^n-U_{\delta}^{n-1},1\right\rangle_h \left\langle 1,U_{\delta}^{n+1}-U_{\delta}^{n}\right\rangle_h\\
     = & -\frac{M}{2} \left\langle U_{\delta}^{n+1} - \omega,1\right\rangle_h^2 +\frac{M}{2} \left\langle U_{\delta}^n -\omega ,1\right\rangle_h ^2\\
        & +\frac{M}{2} \left\langle 1,U_{\delta}^{n+1}-U_{\delta}^{n}\right\rangle_h^2-M\left\langle U_{\delta}^n-U_{\delta}^{n-1},1\right\rangle_h \left\langle 1,U_{\delta}^{n+1}-U_{\delta}^{n}\right\rangle_h.
\end{align*}
Introducing the discrete counterpart of the penalized NOK energy (\ref{functional:pNOK}), 
\begin{align}
E^{\text{pNOK}}_{h}(U^{n}_{\delta}) =  \frac{\epsilon}{2}\left\langle\mathcal{L}_{\delta,h}U^{n}_{\delta},U^{n}_{\delta}\right\rangle_{h}+\frac{1}{\epsilon}\left\langle W(U^{n}_{\delta}),1 \right\rangle_h + \frac{\gamma}{2}\left\langle\mathcal{L}_{\delta,h}^{-1}[U^{n}_{\delta}-\omega],U^{n}_{\delta}-\omega\right\rangle_h + \frac{M}{2} \left\langle U^{n}_{\delta}, 1 \right\rangle_h^2,
\end{align}
and reordering terms from the two sides, then we have
\begin{align*}
    & \left[E_h^{\text{pNOK}}(U_{\delta}^{n+1})+\left(\frac{A_h}{2}+\frac{1}{4\tau}\right)\|U_{\delta}^{n+1}-U_{\delta}^n\|_{L^2,h}^2+\frac{\gamma B_h}{2}\left\langle\mathcal{L}^{-1}_{\delta,h}(U_{\delta}^{n+1}-U_{\delta}^{n}),U_{\delta}^{n+1}-U_{\delta}^n\right\rangle_h\right] \\
    & -\left[E_h^{\text{pNOK}}(U_{\delta}^n)+\left(\frac{A_h}{2}+\frac{1}{4\tau}\right)\|U_{\delta}^{n}-U_{\delta}^{n-1}\|_{L^2,h}^2+\frac{\gamma B_h}{2}\left\langle\mathcal{L}^{-1}_{\delta,h}(U_{\delta}^{n}-U_{\delta}^{n-1}),U_{\delta}^{n}-U_{\delta}^{n-1}\right\rangle_h\right] \\
    & +\frac{1}{\tau}\|U_{\delta}^{n+1}-U_{\delta}^n\|_{L^2,h}^2+\left(\frac{A_h}{2}+\frac{1}{4\tau}\right)\|U_{\delta}^{n+1}-2U_{\delta}^n+U_{\delta}^{n-1}\|_{L^2,h}^2 \\
    & +\frac{\gamma B_h}{2}\left\langle\mathcal{L}^{-1}_{\delta,h}(U_{\delta}^{n+1}-2U_{\delta}^{n}+U_{\delta}^{n-1}),U_{\delta}^{n+1}-2U_{\delta}^{n}+U_{\delta}^{n-1}\right\rangle_h+\frac{\epsilon}{2}\left\langle\mathcal{L}_{\delta,h}(U_{\delta}^{n+1}-U_{\delta}^n),U_{\delta}^{n+1}-U_{\delta}^n\right\rangle_h \\
   = \ &  \frac{W''(\xi_h^n)}{2\epsilon}\|U_{\delta}^{n+1}-U_{\delta}^n\|_{L^2,h}^2-\frac{W''(\mu_h^n)}{\epsilon}\left\langle U_{\delta}^n-U_{\delta}^{n-1},U_{\delta}^{n+1}-U_{\delta}^n\right\rangle_h \nonumber\\
    & +\frac{\gamma}{2}\left\langle\mathcal{L}^{-1}_{\delta,h}[U_{\delta}^{n+1}-U_{\delta}^n],U_{\delta}^{n+1}-U_{\delta}^n\right\rangle_h- \gamma\left\langle\mathcal{L}^{-1}_{\delta,h}[U_{\delta}^n-U_{\delta}^{n-1}],U_{\delta}^{n+1}-U_{\delta}^n\right\rangle_h \nonumber \\
    & +\frac{M}{2}\left|\left\langle U_{\delta}^{n+1}-U_{\delta}^n,1\right\rangle_h\right|^2 -M\left\langle U_{\delta}^{n}-U_{\delta}^{n-1},1\right\rangle_h\left\langle U_{\delta}^{n+1}-U_{\delta}^n,1\right\rangle_h \\
   \leq \ & \frac{L_{W''}}{2\epsilon}\|U_{\delta}^{n+1}-U_{\delta}^n\|_{L^2,h}^2 + \frac{L_{W''}}{2\epsilon}\left(\|U_{\delta}^{n}-U_{\delta}^{n-1}\|_{L^2,h}^2+\|U_{\delta}^{n+1}-U_{\delta}^n\|_{L^2,h}^2\right) \\
    & + \frac{\gamma}{2} \|\mathcal{L}^{-1}_{\delta,h}\|_{L^2,h}\|U_{\delta}^{n+1}-U_{\delta}^n\|_{L^2,h}^2 + \frac{\gamma}{2}\|\mathcal{L}^{-1}_{\delta,h}\|_{L^2,h}\left(\|U_{\delta}^{n}-U_{\delta}^{n-1}\|_{L^2,h}^2+\|U_{\delta}^{n+1}-U_{\delta}^n\|_{L^2,h}^2\right) \\
    & + \frac{M}{2}|\Omega|\|U_{\delta}^{n+1}-U_{\delta}^n\|_{L^2,h}^2 + \frac{M}{2}|\Omega|\left(\|U_{\delta}^{n}-U_{\delta}^{n-1}\|_{L^2,h}^2+\|U_{\delta}^{n+1}-U_{\delta}^n\|_{L^2,h}^2\right) \\
    = \ & C_h \left(\|U_{\delta}^{n}-U_{\delta}^{n-1}\|_{L^2,h}^2+2\|U_{\delta}^{n+1}-U_{\delta}^n\|_{L^2,h}^2\right)
\end{align*}
where the constant $C$ is given by 
\begin{align*}
    C_h=\frac{L_{W''}}{2\epsilon}+\frac{\gamma}{2}\|\mathcal{L}^{-1}_{\delta,h}\|_{L^2,h}+\frac{M}{2}|\Omega|.
\end{align*}
Note that $\|\mathcal{L}^{-1}_{\delta,h}\|_{L^2,h}$ denote the $L^2$-bound of $\mathcal{L}^{-1}_{\delta,h}$, which is provided in Corollary (\ref{remark:L2bdd}). Adding $C_h(\|U_{\delta}^{n+1}-U_{\delta}^{n}\|_{L^2,h}^2-\|U_{\delta}^{n}-U_{\delta}^{n-1}\|_{L^2,h}^2)$ to both side of the above inequality, we further have
\begin{align*}
    \text{LHS}+C_h\left(\|U_{\delta}^{n+1}-U_{\delta}^{n}\|_{L^2,h}^2-\|U_{\delta}^{n}-U_{\delta}^{n-1}\|_{L^2,h}^2\right)\leq \text{RHS} +  3C_h\|U_{\delta}^{n+1}-U_{\delta}^{n}\|_{L^2,h}^2.
\end{align*}
Subsequently, 
\begin{align}\label{ieq:energystable1}
    & \left[E_h(U_{\delta}^{n+1})+\left(\frac{A_h}{2}+\frac{1}{4\tau}+C_h\right)\|U_{\delta}^{n+1}-U_{\delta}^n\|_{L^2,h}^2+\frac{\gamma B_h}{2}\langle\mathcal{L}^{-1}_{\delta,h}(U_{\delta}^{n+1}-U_{\delta}^{n}),U_{\delta}^{n+1}-U_{\delta}^n\rangle_h\right] \nonumber\\
    & -\left[E_h(U_{\delta}^n)+\left(\frac{A_h}{2}+\frac{1}{4\tau}+ C_h\right)\|U_{\delta}^{n}-U_{\delta}^{n-1}\|_{L^2,h}^2+\frac{\gamma B_h}{2}\langle\mathcal{L}^{-1}_{\delta,h}(U_{\delta}^{n}-U_{\delta}^{n-1}),U_{\delta}^{n}-U_{\delta}^{n-1}\rangle_h\right] \nonumber\\
    & +\frac{1}{\tau}\|U_{\delta}^{n+1}-U_{\delta}^n\|_{L^2,h}^2+\left(\frac{A_h}{2}+\frac{1}{4\tau}\right)\|U_{\delta}^{n+1}-2U_{\delta}^n+U_{\delta}^{n-1}\|_{L^2,h}^2 \nonumber\\
    & +\frac{\gamma B_h}{2}\left\langle\mathcal{L}^{-1}_{\delta,h}(U_{\delta}^{n+1}-2U_{\delta}^{n}+U_{\delta}^{n-1}),U_{\delta}^{n+1}-2U_{\delta}^{n}+U_{\delta}^{n-1}\right\rangle_h+\frac{\epsilon}{2}\left\langle\mathcal{L}_{\delta,h}(U_{\delta}^{n+1}-U_{\delta}^n),U_{\delta}^{n+1}-U_{\delta}^n\right\rangle_h \nonumber\\
    \leq &\ 3C_h\|U_{\delta}^{n+1}-U_{\delta}^{n}\|_{L^2,h}^2.
\end{align}

Let us define 
\begin{align}\label{eqn:modifiedenergy}
\tilde{E}^{\text{pNOK}}_h (U_{\delta}^{n},U_{\delta}^{n-1} ) = E_h^{\text{pNOK}}(U_{\delta}^n)+\left(\frac{A_h}{2}+\frac{1}{4\tau}+C_h\right)\|U_{\delta}^{n}-U_{\delta}^{n-1}\|_{L^2}^2+\frac{\gamma B_h}{2}\|\mathcal{L}^{-\frac{1}{2}}_{\delta,h}\left(U_{\delta}^{n}-U_{\delta}^{n-1}\right)\|_{L^2,h}.
\end{align}
By observing the two sides of (\ref{ieq:energystable1}), we note that when $\frac{A_h}{2}+\frac{1}{4\tau}\ge0$, $B_h\ge0$, and $\frac{1}{\tau}\ge 3 C_h$, it yields the discrete energy stability: $\tilde{E}^{\text{pNOK}}_h (U_{\delta}^{n+1},U_{\delta}^{n} ) \le \tilde{E}^{\text{pNOK}}_h (U_{\delta}^{n},U_{\delta}^{n-1} ) $.

We summarize the above discussion as a theorem.
\begin{theorem}\label{thm:BDF_NOK_energy_fully}

For the second-order fully-discrete BDF scheme in (\ref{eqn:BDF_NOK_fully}), define a modified energy functional as in (\ref{eqn:modifiedenergy}) with $C_h$ given by
\begin{align*}
    C_h=\frac{L_{W''}}{2\epsilon}+\frac{\gamma}{2}\|\mathcal{L}^{-1}_{\delta,h}\|_{L^2,h}+\frac{M}{2}|\Omega|.
\end{align*}
Then we have the energy stability
\[
\tilde{E}^{\emph{pNOK}}_h (U_{\delta}^{n+1},U_{\delta}^{n} ) \le \tilde{E}^{\emph{pNOK}}_h (U_{\delta}^{n},U_{\delta}^{n-1} ), 
\]
provided that $A_h \ge 0, B_h \ge 0$ and $\tau \le \frac{1}{3C_h}$. In particular, if the stabilizers $A_h = B_h = 0$, we have the energy stability in a simpler form:
\[
E_h^{\emph{pNOK}}(U_{\delta}^{n+1})+\left( \frac{1}{4\tau}+C_h\right)\|U_{\delta}^{n+1}-U_{\delta}^{n}\|_{L^2}^2
\le
E_h^{\emph{pNOK}}(U_{\delta}^n)+\left( \frac{1}{4\tau}+C_h\right)\|U_{\delta}^{n}-U_{\delta}^{n-1}\|_{L^2}^2.
\]

\end{theorem}

\begin{remark}
Though there is no influence on the discrete energy stability whether or not we incorporate the stabilizers $A_h$ and $B_h$ in the 2nd-order BDF scheme, we find in the numerical implementation that adding two stabilizers in the scheme will allow us to take a much larger time step $\tau$ in practice. Therefore in Section \ref{sec:es}, we will pick $A_h>0$ and $B_h>0$ for all the numerical simulations.
\end{remark}

\section{Numerical Experiments}\label{sec:ne}

In this section, we use the second-order BDF scheme (\ref{eqn:BDF_NOK_fully}) to solve the pNACOK equation (\ref{eqn:pNOK}) with periodic boundary conditions. In our numerical experiments, we apply the ODE solver proposed in \cite{Du_Yang2017} to numerically compute the Fourier symbols $\{\lambda_{\delta}(\mathbf{k})\}$ for the nonlocal operator $\mathcal{L}_{\delta}$. We consider two special kernel functions for the nonlocal operator:
\begin{align}
&\text{Power kernel}: \quad \rho_{\delta}(s) = \frac{C_{\alpha}}{\delta^{d+2-\alpha}|s|^{\alpha}}, \ s \in [-\delta,0)\cup(0,\delta], \quad \text{where\ } \alpha \in [0,d+2), \label{eqn:PowerKernel}\\
&\text{Gaussian kernel}: \quad \rho_{\delta}(s) = \frac{4}{\pi^{d/2}\delta^{d+2}} e^{-\frac{s^2}{\delta^2}}, \label{eqn:GaussianKernel}
 \end{align}
in which $C_{\alpha}$ is the normalization constant for the power kernel.

In the numerical experiments, we mainly focus on the 2D examples as they are more relevant to the diblock copolymer system. Therefore, we fix the domain $\Omega=[-\pi, \pi]^2\in\mathbb{R}^2$. Unless  otherwise specified, we fix $N=N_x=N_y= 512$, the mesh size is given by $h=h_x=h_y=\frac{\pi}{256}$, $\epsilon = 10 h$, step size $\tau = 10^{-3}$ and $M=1000$. Other parameters, such as $\omega$, $\gamma$, $\alpha$, $\delta$, $A_h$ and $B_h$ vary in different simulations. The stopping criteria for the time iteration are set to be
\begin{align}\label{eqn:stop_condition}
    \frac{||U_{\delta}^{n+1}-U_{\delta}^{n}||_{h,L^{\infty}}}{\tau} \leq 10^{-5}.
\end{align}

\subsection{Asymptotic Compatibility}

In this subsection, we numerically validate the asymptotic compatibility for the pNACOK equation.  To this end, we take a small disk centered at origin as the initial data:
\begin{align*}
        u^{(0)}(x,y) = \left\{
    \begin{array}{ll}
       1  &  \text{if } x^2 + y^2 < 4\omega/\pi,\\
       0  &  \text{otherwise},
    \end{array}
    \right. 
\end{align*}
We take $\omega = 0.1$, $\gamma = 100$, $A_h = 5000$ and $B_h = 5$. The nonlocal operator $\cL_{\delta}$ is chosen to be with the power kernel $\alpha = 0$ and nonlocal horizon $\delta = 0.5$. We numerically test the errors between nonlocal and local models by taking a small time step size $\tau = 10^{-3}$ and various values of $\delta$. We simulate the pNACOK equation up to time $T=0.01, 1$ and $10$.

Table \ref{table:NACOK_AC} presents the errors and corresponding convergence rate at a fixed time $T=0.01, 1$ and $10$. It can be observed that the numerical convergence rate is approximately $O(\delta^2)$,  which is consistent with the theoretical prediction of asymptotic compatibility of the  pNACOK equation.

\begin{table}[h!]
\begin{center}
\begin{tabular}{ |c||c|c|c|c|c|c|  }
\hline
 & \multicolumn{2}{|c|}{$T = 0.01$} & \multicolumn{2}{|c|}{$T = 1$} & \multicolumn{2}{|c|}{$T = 10$} \\
\hline
$\delta = 1$ & $\|u_h^{\delta} - u_h^0\|_2$ & Rate & $\|u_h^{\delta} - u_h^0\|_2$ & Rate & $\|u_h^{\delta} - u_h^0\|_2$ & Rate \\
\hline
$\delta$ & 3.6365e-3 & -       & 9.9332e-1 &  - & 1.0534 &  -\\
$\delta/2$ & 8.8995e-4 & 2.0268 & 2.0916e-1 & 2.2476 & 2.2437e-1 &  2.2311\\
$\delta/4$ & 2.1714e-4 & 2.0351 & 4.4585e-2 & 2.2300 & 4.8252e-2 &  2.2172\\
$\delta/8$ & 5.1960e-5 & 2.0633 & 8.5567e-3 & 2.3814 & 9.2525e-3 &  2.3827\\
\hline
\end{tabular}
\end{center}
\caption{Errors between numerical solutions of the pNACOK and pLACOK equations at time $T=0.01,1$ and $10$. Other parameters are $\omega = 0.1$, $\gamma = 100$, $\alpha=0$, $\delta=0.5$, $A_h = 5000$ and $B_h = 5$.}
\label{table:NACOK_AC}
\end{table}

\begin{figure}[htbp]
\centerline{
\includegraphics[width=0.6\textwidth]{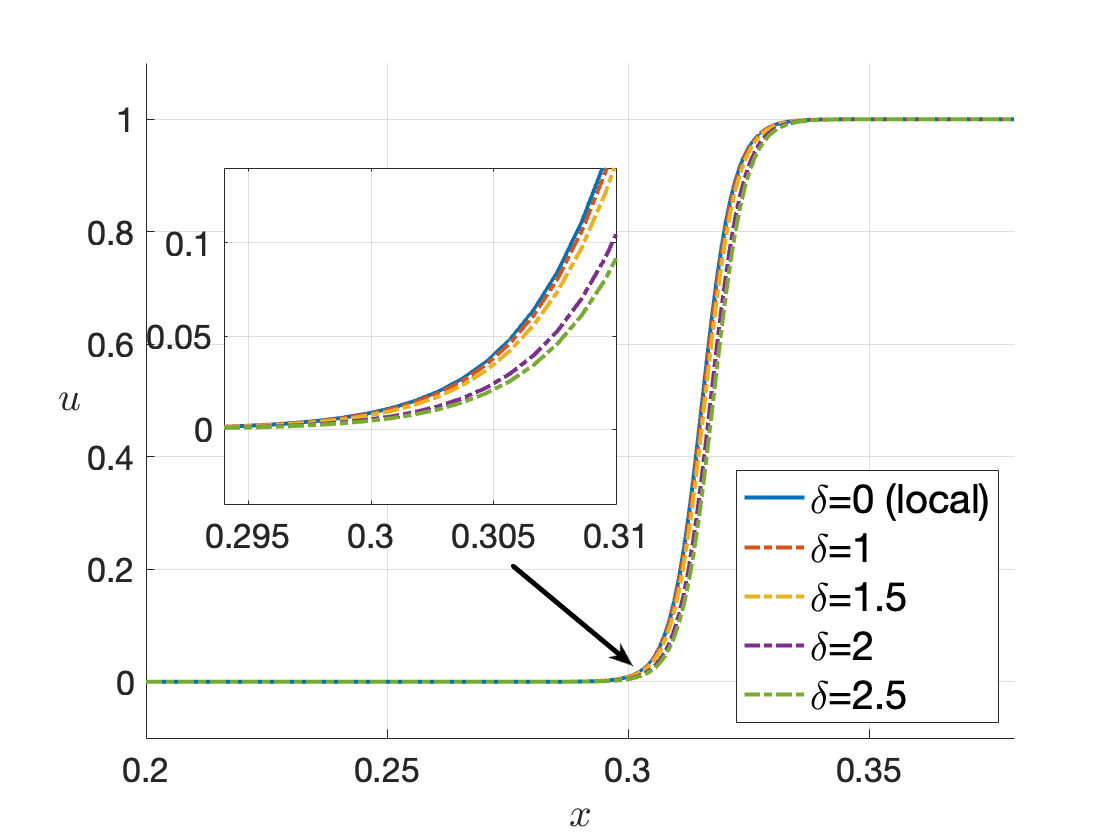}
$ \quad $ }
\vspace{-4 mm}
\caption{Comparison between the numerical solutions of pNACOK and pLACOK equations for vary values of nonlocal horizon $\delta=1.0, 1.5, 2.0, 2.5$. Other parameters are $\omega = 0.3$, $\gamma = 200$, $\alpha=0$, $\delta=0.5$, $A_h = 2000$ and $B_h = 0$.}
\label{fig:ac_dynamics}
\end{figure}

In Figure (\ref{fig:ac_dynamics}),  we present a numerical example to compare the solutions of the pNACOK and pLACOK equations in 1D. In this example, the initial data is $u^{(0)}(x) = 0.2\sin(\frac{\pi x}{4})$, time step is $\Delta t = 10^{-3}$. The nonlocal operator is with power-law kernel and $\alpha = 0$, $\delta = 0.5$. We take a finer mesh resolution $N = 2048$, and solve the pNACOK and pLACOK equations with various value $\delta = 1.0, 1.5, 2.0, 2.5$ up to time $T = 100$. Figure (\ref{fig:ac_dynamics}) shows that as $\delta$ becomes smaller, the numerical solution of the pNACOK equation becomes closer to that of the pLACOK equation.

%

\subsection{Temporal Convergence Rate}
Now we numerically verify the temporal rate of convergence for the 2nd-order BDF scheme (\ref{eqn:BDF_NOK_fully}). To this end, we consider the 1D pNACOK equation and take the initial data 
\begin{align*}
        u_0(x) = \left\{
    \begin{array}{ll}
       1  &  \text{if } (x-0.1)^2  < \omega/(2\pi)+0.03,\\
       0  &  \text{otherwise},
    \end{array}
    \right. 
\end{align*}
and run the numerical simulation up to time $T = 0.02$. we fix $\omega=0.3$, $\gamma=200$, $A_h = 2000$, $B_h = 0$. Taking the nonlocal operator $\mathcal{L}_{\delta}$ be with the power kernel, we fix $\alpha=2$ and $\delta=0.5$. We generate a benchmark solution by the 2nd-order BDF scheme (\ref{eqn:BDF_NOK_fully}) with a tiny time step $\tau = 1e-6$. then we compute the discrete $L^2$ error between the numerical solutions with larger step sizes and the benchmark one. 
    
Table \ref{table:Time_Conv} shows the errors and the convergence rates of second-order BDF scheme (\ref{eqn:BDF_NOK_fully}) based on the data at $T=0.02$ with the time steps $\tau = \frac{10^{-3}}{1}, \frac{10^{-3}}{2}, \frac{10^{-3}}{4}, \cdots, \frac{10^{-3}}{32}$. We verify the convergence rates for two different values of $\epsilon = 10h$ and $20h$ with $h = \frac{2\pi}{N}$ and $N = 512$. From the table, we can observe that the numerical rates tend to approach the theoretical value 2. 

\begin{table}[h!]
\begin{center}
\begin{tabular}{ |c||c|c|c|c|  }
\hline
 & \multicolumn{2}{|c|}{$\epsilon =10h$} & \multicolumn{2}{|c|}{$\epsilon =20h$}  \\
\hline
$\Delta t = $1e-3 & Error & Rate & Error & Rate  \\
\hline
$\Delta t $ & 1.157598 & -       & 0.672051 &  - \\
$\Delta t/2 $ & 0.164550 & 2.814531 & 0.325463 & 1.046075 \\
$\Delta t/4 $ & 0.106904 & 0.622206 & 0.066833 & 2.283852 \\
$\Delta t/8 $ & 0.032665 & 1.710477& 0.013439 & 2.314069 \\
$\Delta t/16 $ & 0.008102 & 2.011280 & 0.003056 & 2.136332 \\
$\Delta t/32 $ & 0.001811 & 2.161469 & 0.000647 & 2.238881 \\
1e-6 (Benchmark)  & - & - & - & - \\
\hline
\end{tabular}
\end{center}
\caption{The convergence rate for the 2nd-order BDF scheme (\ref{eqn:BDF_NOK_fully}) at time $T = 0.02$. Other parameters are $\omega = 0.3$, $\gamma = 200$, $\alpha=2$, $\delta=0.5$, $A_h = 2000$, and $B_h = 0$.}
\label{table:Time_Conv}
\end{table}

\begin{figure}[htbp]
\begin{center}
\includegraphics[width=0.7\textwidth]{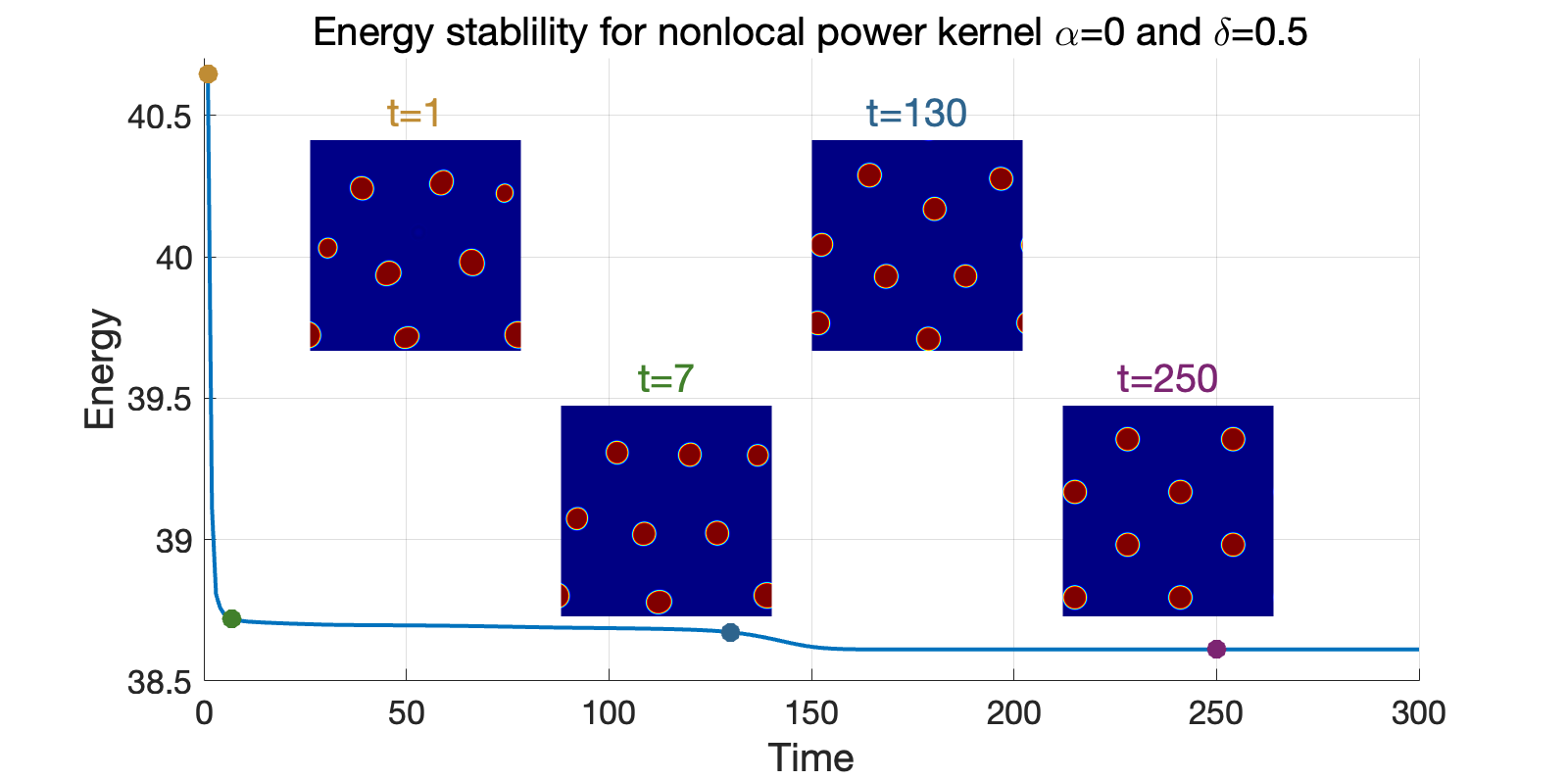}
\includegraphics[width=0.7\textwidth]{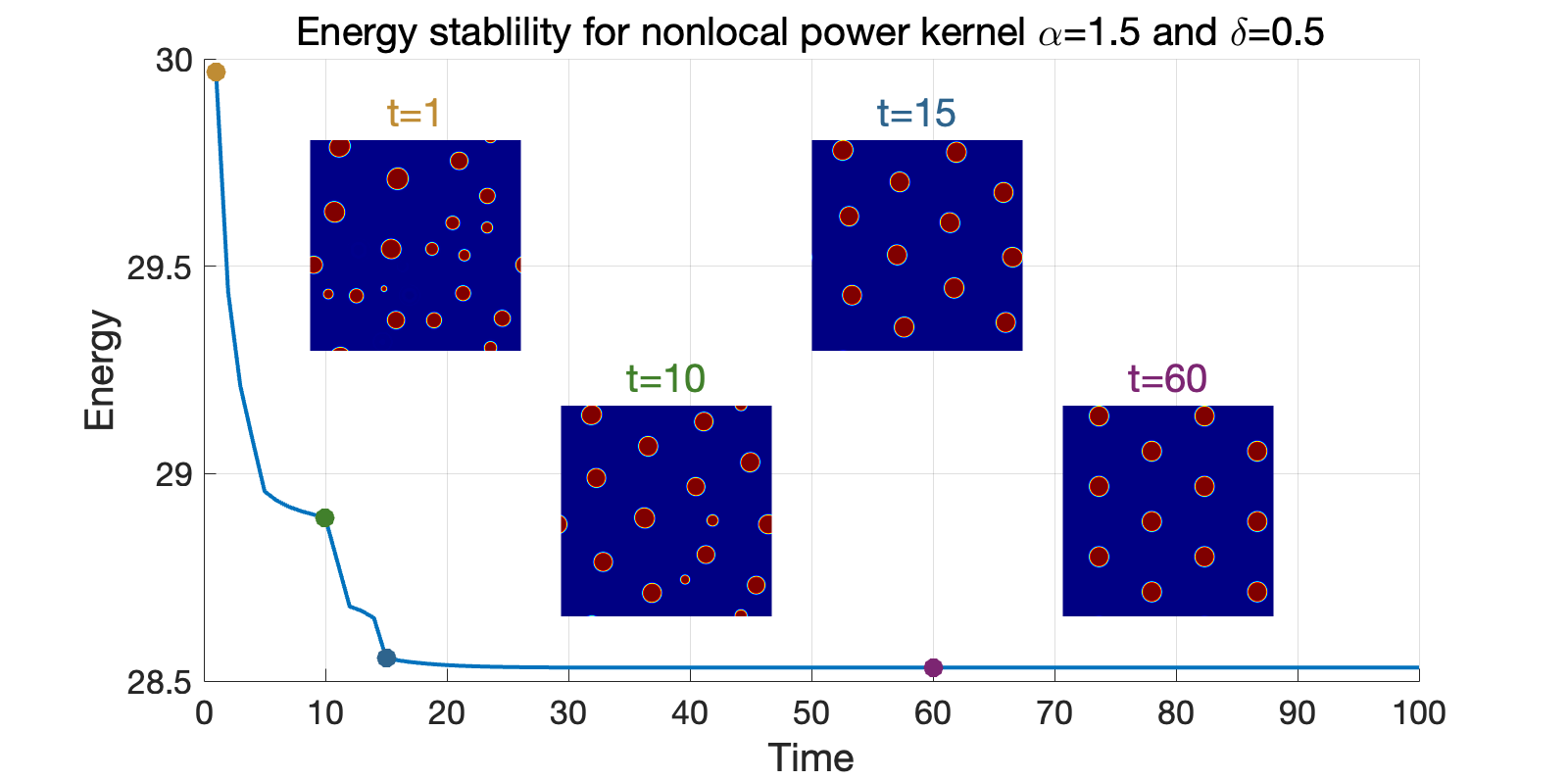}
\end{center}
\caption{Two coarsening dynamic processes show the energy stability for the BDF schemes (\ref{eqn:BDF_NOK_fully}) and different patterns at equilibria. (Top). Square lattice is formed at equilibrium. Here $\alpha = 0$; (Bottom). Hexagonal lattice is formed at equilibrium. Here $\alpha = 1.5$. Other parameter values are $\delta = 0.5$, $\omega = 0.1$, $\gamma = 7000$, $A_h=5000$, and $B_h=5$. }
\label{fig:energy_stability}
\end{figure}

\subsection{Coarsening Dynamics and Energy Stability}

In this subsection, we focus on the coarsening dynamics and the energy stability of the pNACOK equation implementing by the BDF scheme (\ref{eqn:BDF_NOK_fully}). 

In our numerical experiments, we take a random initial which is generated on a uniform mesh over the domain $\Omega$ with mesh size being $16h$. This can be achieved in Matlab by the command 
\[
\texttt{repelem}(\texttt{rand}(\texttt{N}/\texttt{ratio},\texttt{N}/\texttt{ratio}),\texttt{ratio},\texttt{ratio})
\]
with $\texttt{ratio} = 16$. Meanwhile, we choose the time step $\tau = 10^{-3}$ and terminate the implementation by the stopping condition (\ref{eqn:stop_condition}). For the other parameters, we fix $\omega = 0.1$, $\gamma = 7000$, $A_h = 5000$, and $B_h = 5$.

We simulate the pNACOK equation for two different nonlocal operators. One is with the power-law kernel and $\alpha = 0$, the other is with the power-law kernel and $\alpha = 1.5$. We fix the nonlocal horizon $\delta = 0.5$.

The coarsening dynamics and the energy stability for the pNACOK equation and its equilibria are shown in Figure (\ref{fig:energy_stability}). For each subfigure, starting from the random initial, we insert the snapshots taken at four different times $t$, with colored titles corresponding to the colored marker on the monotone decreasing energy curve. The phase separations arise in a very short time period and a group of bubbles with different sizes appear, then the small bubbles disappear, other bubbles evolves into being of equal size. Eventually all the equally-sized bubbles become equally distanced, and form different bubble structures in the 2D domain $\Omega$. In the top subfigure, with $\alpha = 0$, the equilibria becomes square lattice; while in the bottom subfigure, with $\alpha = 1.5$, it forms a hexagonal pattern. A more detailed study on the different patterns of the equilibria for the nonlocal OK model is systematically explored in our ongoing work \cite{Zhao_Luo2d}. In addition, the energy curve clearly shows the energy stability of the BDF scheme (\ref{eqn:BDF_NOK_fully}) for the pNACOK equation.

\begin{figure}[htbp]
\centerline{
\includegraphics[width=0.7\textwidth]{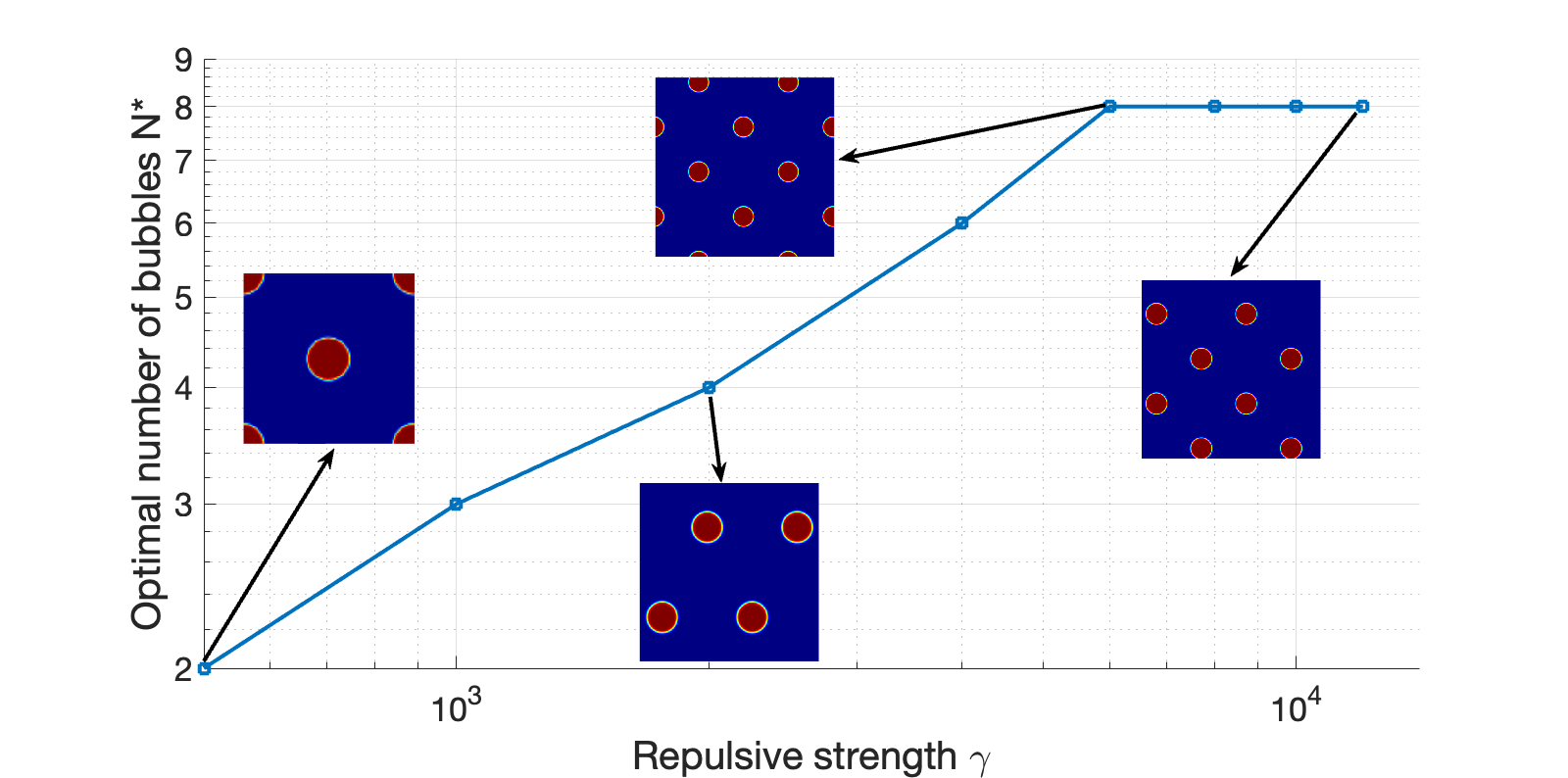}
$ \quad $ }
\vspace{-0 mm}
\caption{The effect of $\gamma$ on the optimal number $N^*$ of bubbles for the NOK system. Here we take the nonlocal operator to be with the power-law kernel, and $\alpha=0$, $\delta=0.5$. For various $\gamma = 500, 1000, 2000, 4000, 6000, 8000,10000, 12000$, the corresponding optimal number of bubbles are $N^* = 2, 3, 4, 6, 8, 8, 8, 8$, respectively. The four insets are equilibrium states for $\gamma = 500, 2000, 6000,12000$. Other parameter values are $\omega=0.1$, $A_h = 5000$, and $B_h = 5$.}
\label{fig:gamma_increase}
\end{figure}

\begin{figure}[htbp]
\centerline{
\includegraphics[width=0.5\textwidth]{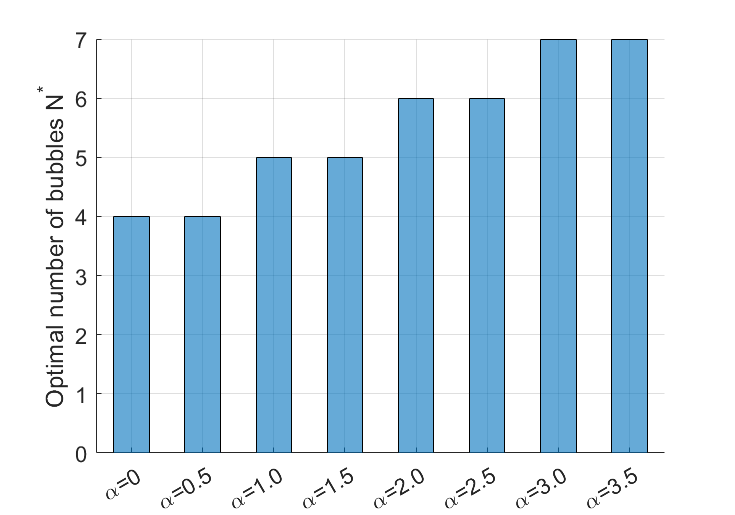}
$ \quad $ }
\vspace{-0 mm}
\caption{The $\alpha$-effect on the optimal number of bubbles $N^*$. Here the nonlocal operator is with the power-law kernel. For various $\alpha = 0, 0.5, 1.0, 1.5, 2.0, 2.5, 3.0, 3.5$, the optimal number of bubbles at equilibria is $N^* = 4, 4, 5, 5, 6, 6, 7, 7$, respectively. We fix $\gamma = 2000$. Other parameter values are the same as the previous example.}
\label{fig:optimal}
\end{figure}

\subsection{The $\alpha$-effect on the Equilibria for the Nonlocal Operator with Power-law Kernels}

Now we study the system equilibria for the nonlocal operator with power-law kernels, but with various long-range interaction strength $\gamma$. Different from the original OK model in \cite{OhtaKawasaki_Macromolecules1986}, in which the number of bubbles increase as $\gamma \rightarrow \infty$, the nonlocal OK model gives a upper bound for the number of bubbles even when $\gamma \rightarrow \infty$.
 
 Here we take power-law kernel with $\alpha=0$ and $\delta=0.5$. We fix the parameters $\omega=0.1$, $A_h=5000$, $B_h=5$, and various $\gamma$ values from $500$ to $12000$. Figure \ref{fig:gamma_increase} shows that as $\gamma$ increases, the optimal number $N^*(\gamma)$ of bubbles increases. When $\gamma$ exceeds certain threshold, $N^*(\gamma)$ stop growing and stay as a constant.

Figure \ref{fig:gamma_increase} shows that there exists an upper bound $\bar{N}$ for the optimal number of bubbles at equilibria when $\gamma\rightarrow$. In the next numerical experiment, we study the influence of $\alpha$ on the optimal number $N^*$. Here we use the same parameter values for $\epsilon, A_h, B_h, \omega$ and $\delta$, but consider the power-law kernels with various $\alpha = 0, 0.5, 1.0, 1.5, 2.0, 2.5, 3.0, 3.5$. For a fixed $\gamma = 2000$, Figure \ref{fig:optimal} shows that the optimal number of bubbles increases as $\alpha$ increases.

\begin{figure}[htbp]
\centerline{
\includegraphics[width=0.7\textwidth]{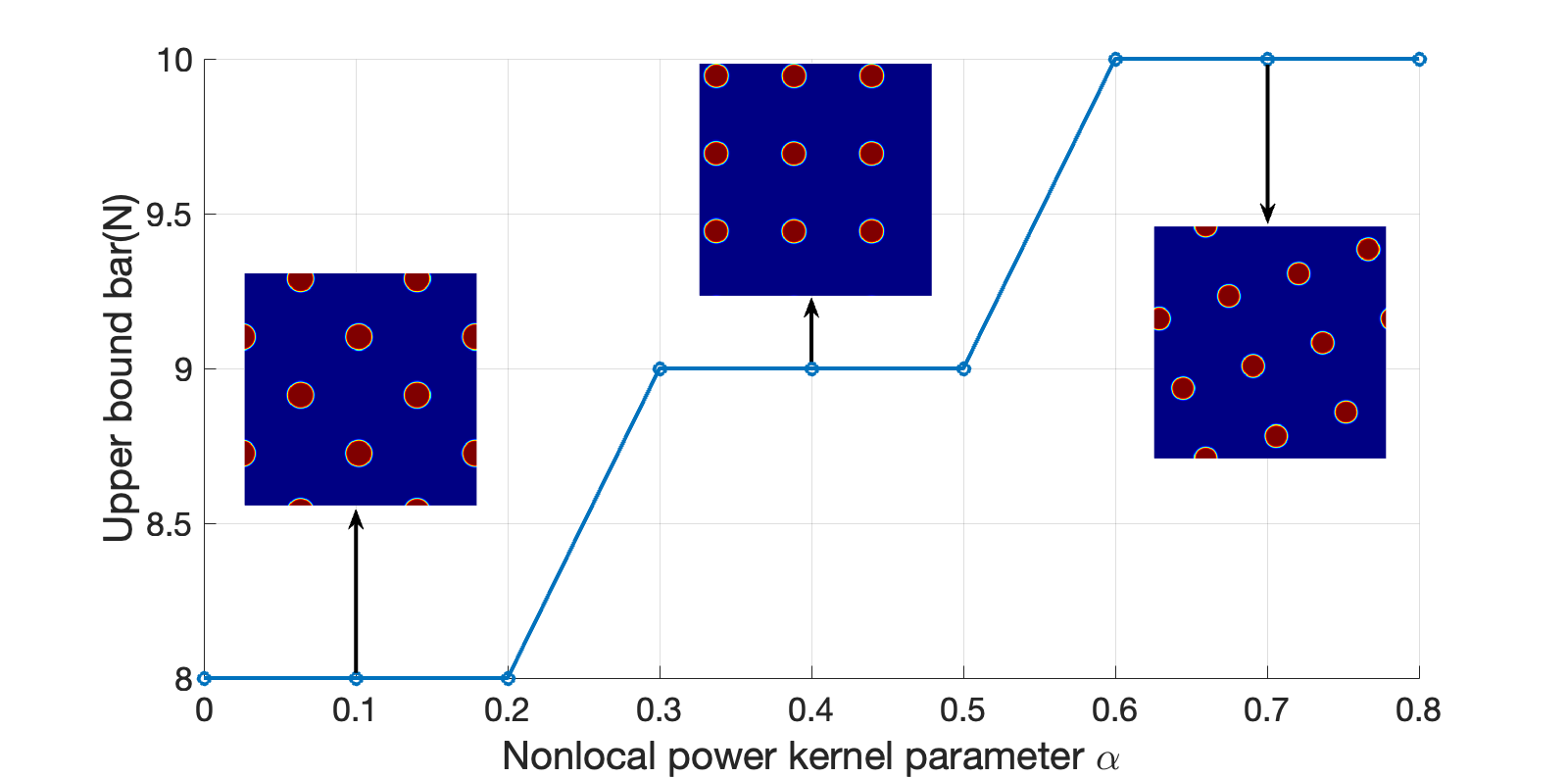}
$ \quad $ }
\vspace{-0 mm}
\caption{The $\alpha$-effect on the upper bound $\bar{N}$ of the optimal number of bubbles at equilibria. For $\alpha = 0, 0.1, 0.2, 0.3, 0.4, 0.5, 0.6, 0.7, 0.8$ and corresponding upper bound $\bar{N} = 8,  8,  8,  9,  9,  9,  10,  10, 10$, respectively. We fix $\omega = 0.1$ and $\delta=0.5$. For each simulation, we take sufficient large $\gamma$ so that the optimal number of bubbles reaches the upper bound. Other parameters remain the same as the previous example.}
\label{fig:upper1}
\end{figure}

\begin{figure}[htbp]
\centerline{
\includegraphics[width=0.7\textwidth]{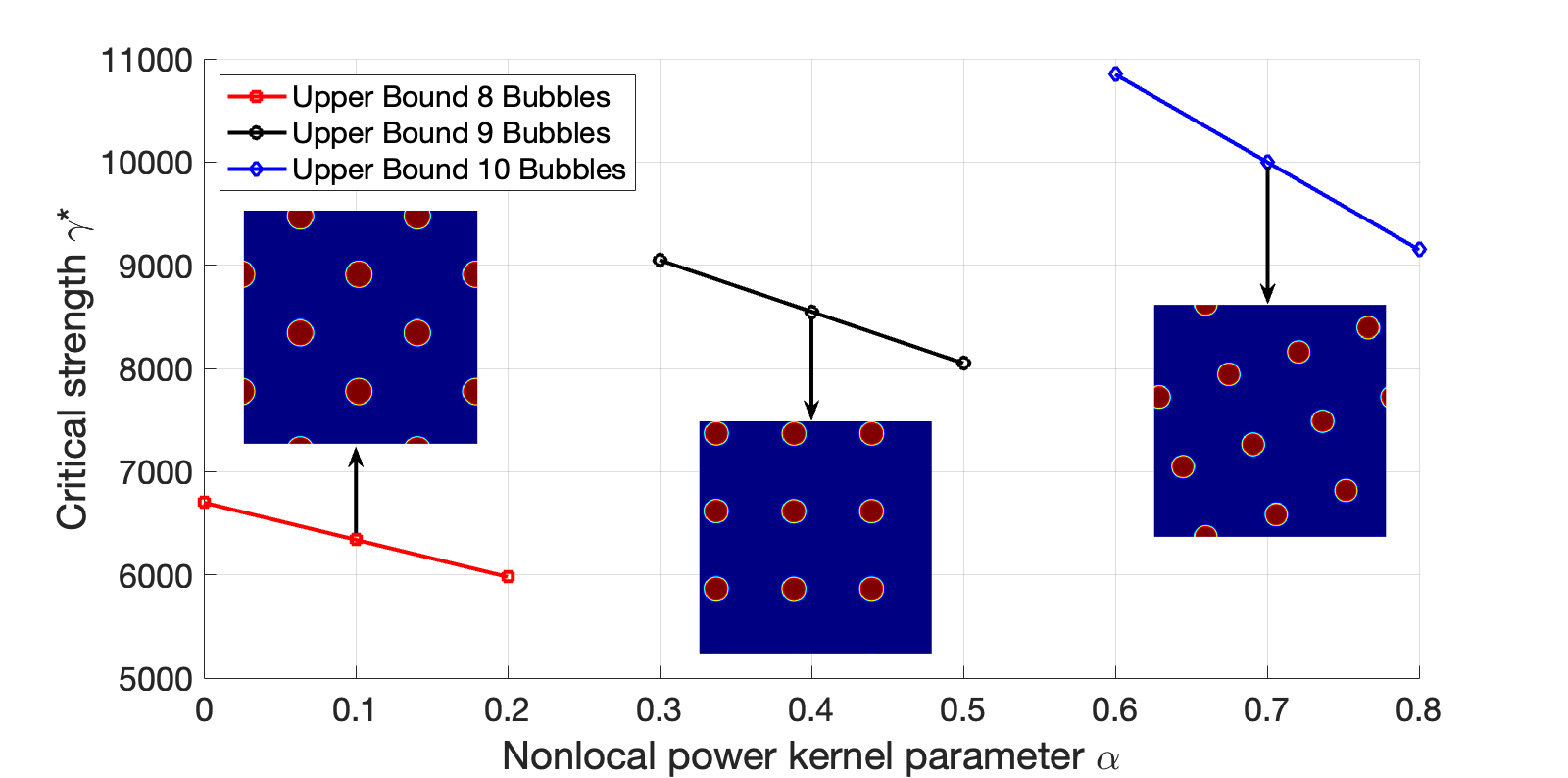}
$ \quad $ }
\vspace{-0 mm}
\caption{The $\alpha$-effect on the critical value of the long-range interaction strength $\gamma^*$.}
\label{fig:upper2}
\end{figure}

Figure \ref{fig:optimal} shows that for a fixed $\gamma$, the optimal number of bubbles $N^*$ increases as $\alpha$ increases. In the next numerical experiment, we further study how $\alpha$ influences the upper bound $\bar{N}$ for the optimal number of bubbles. For each value of $\alpha$, we also find the critical value $\gamma^*$, beyond which the optimal number of bubbles will remain $\bar{N}$. This result is shown in Figure \ref{fig:upper1} and Figure \ref{fig:upper2}. In Figure ref{fig:upper1}, we see that when $\alpha$ increases, the upper bound $\bar{N}$ is nondecreasing. In Figure \ref{fig:upper2}, we still take the same set of values for $\alpha$, but calculate the critical value $\gamma^*$. Interestingly, we find that though the upper bound $\bar{N}$ is nondecreasing, but in the range of $\alpha$ when $\bar{N}$ is unchanged (say, for instance $0\le \alpha \le 0.2$), the critical value $\gamma^*$ decreases.

\begin{figure}[htbp]
\begin{center}
\includegraphics[width=0.27\textwidth,height=0.23\textwidth]{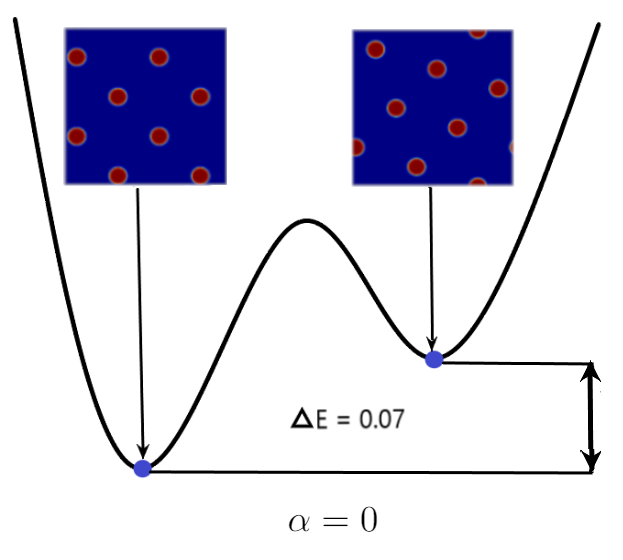}
\includegraphics[width=0.27\textwidth,height=0.23\textwidth]{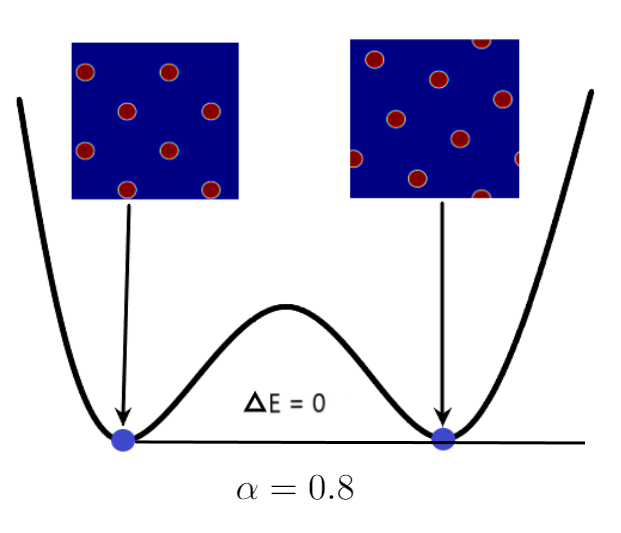}
\includegraphics[width=0.27\textwidth,height=0.23\textwidth]{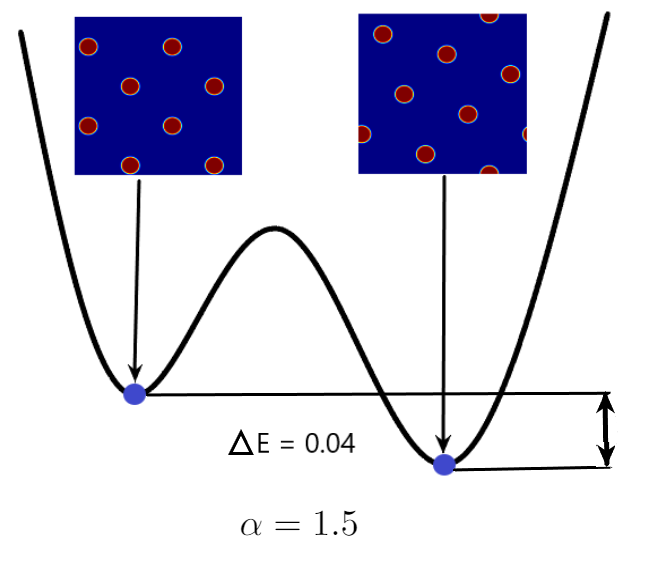}
\end{center}
\caption{The $\alpha$-effect on the pattern: square lattice or hexagonal lattice? }
\label{fig:phase_trans}
\end{figure}

In the previous figures, we observe that for some values of $\alpha$, the equilibria becomes a square lattice pattern; while for other values of $\alpha$, the system displays a hexagonal lattice pattern. In the next example, we numerically verify that for various values of $\alpha$, both square and hexagonal lattice patterns are (local) minima. When $\alpha$ is small, the system prefers a square lattice as the global minimum; while when $\alpha$ is relatively large, the system favors the hexagonal lattice as the global minimum. Figure \ref{fig:phase_trans} present this result. When $\alpha = 0$, the system favors square lattice over hexagonal lattice. When $\alpha = 1.5$, it is the other way around. When $\alpha = 0.8$, both square lattice and hexagonal lattice are global minima.

\subsection{The $\delta$-effect on the Bubble Assembly at Equilibria}

\begin{figure}[htbp]
\begin{center}
\includegraphics[width=0.7\textwidth]{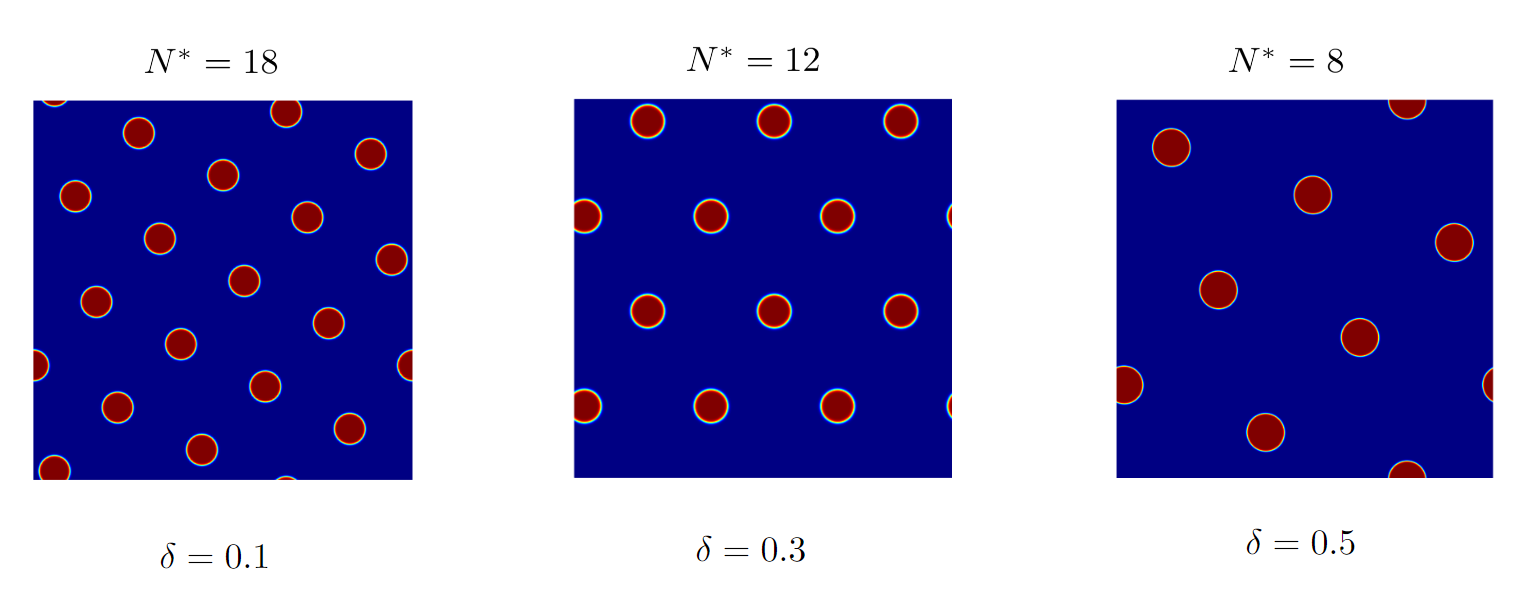}
\includegraphics[width=0.7\textwidth]{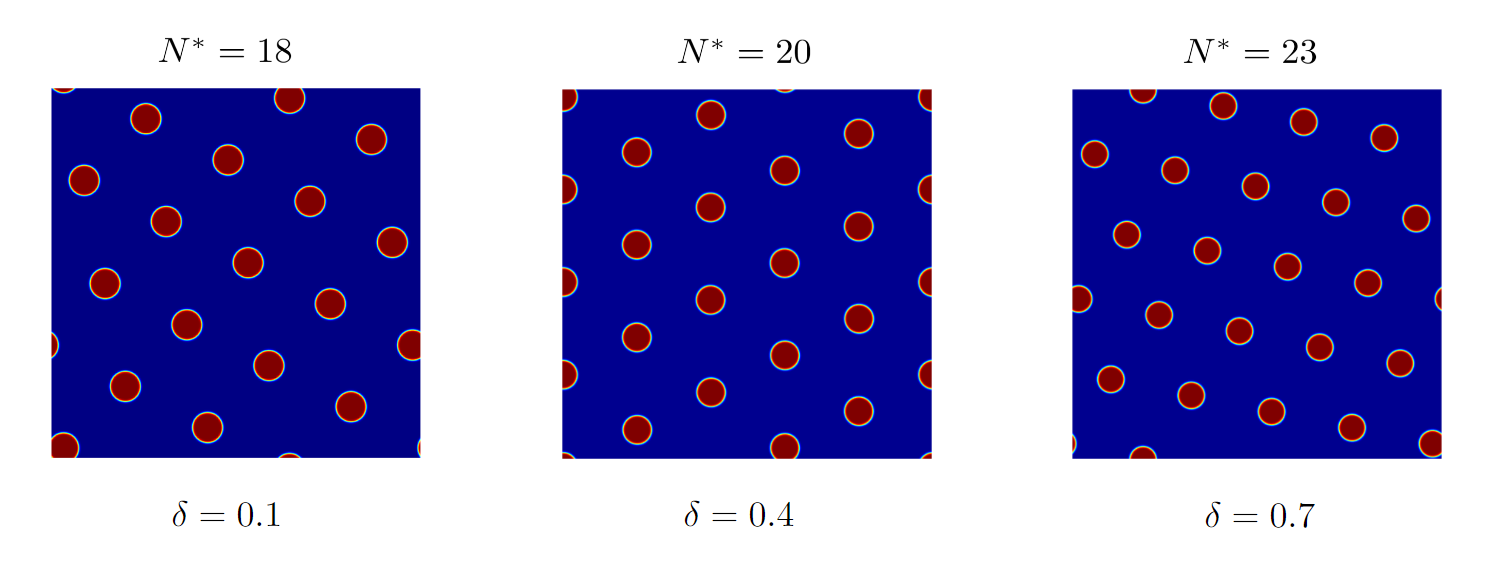}
\end{center}
\caption{Top: the nonlocal horizon $\delta$ demotes the bubble splitting for the Gaussian kernel case Bottom: the nonlocal horizon $\delta$ promotes the bubble splitting for the Power-law kernel case.}
\label{fig:gaussian}
\end{figure}

In this subsection, we explore the influence of the nonlocal horizon $\delta$ on the bubble pattern at equilibria. Here are choose two nonlocal operators. One is with the power-law kernel and $\alpha = 3.5$, the other is with the Gaussian kernel. We fix $\omega = 0.1$, $\gamma = 5000$, $A_h = 5000$, and $B_h = 5$. In the top row of Figure \ref{fig:gaussian}, we increase the value of $\delta = 0.1, 0.3, 0.5$, and observe that the optimal number of bubbles $N^* = 18, 12, 8$, respectively, in a decreasing order. This phenomenon implies that the nonlocal horizon $\delta$ for the Gaussian kernel demotes the bubble splitting. On the other hand, in the bottom row of Figure \ref{fig:gaussian}, as we increase the value of $\delta = 0.1, 0.4, 0.7$, the optimal number of bubbles $N^*$ increases, indicating that $\delta$ promotes the bubble splitting in the power-law kernel case when $\alpha = 3.5$. These numerical findings of the $\delta$-effect on the promotion/demotion of the bubble splitting are consistent with the theoretical results in our previous work \cite{Zhao_Luo1d}.

\section{Concluding remarks}\label{sec:conrem}

In this work, we study the asymptotic compatibility of the spectral collocation methods for the penalized nonlocal ACOK equation in 2D and 3D and the energy stability of the 2nd-order fully-discrete BDF scheme (\ref{eqn:BDF_NOK_fully}).

Meanwhile, we numerically justify the asymptotic compatibility, the convergence rate, and the energy stability of the pNACOK equation by using the proposed scheme (\ref{eqn:BDF_NOK_fully}). Meanwhile, by taking the nonlocal operator as the one with power-law kernel, we examine the $\alpha$-effect on the pattern of bubble assemblies at equilibria for the NOK system. The numerical results indicates that for some values of $\alpha$, the NOK displays a square lattice pattern; for other values of $\alpha$, it favors a hexagonal lattice pattern. Meanwhile, for some $\alpha$, the bubble assembly for the NOK system may have an upper bound on the optimal number of bubbles. These new numerical findings are consistent with our theoretical study in \cite{Zhao_Luo1d}.

There are several directions along which we can extend our work in the future. One possible direction is to perform a systematic study on the pattern formation of the NOK model in 2D. We can examine how the parameter $\delta$ and $\alpha$ affect the pattern formation, either a square lattice or a hexagonal lattice. Another direction is to extend the NOK model to the nonlocal Nakazawa-Ohta (NNO) model. This model is worth comprehensive theoretical and numerical studies. Theoretically, we can conduct 1D analysis on various patterns and understand how the system parameters influence the patterns. Numerically, we can simulate the NNO system and find some possible novel patterns on either square domain or on disk domain. Besides, OK/NOK/NO/NNO models on spherical domain are another important extension, which may be relevant to the deformation of the biomembrane in math biology.


\section{Acknowledgements}

W. Luo's work is supported by NSF grant DMS-2142500. Y. Zhao's work is supported by a grant from the Simons Foundation through Grant No. 357963 and NSF grant DMS-2142500.




\begin{thebibliography}{99}
\itemsep=0pt

\bibitem{OhtaKawasaki_Macromolecules1986} T. Ohta and K. Kawasaki, \emph{Equilibrium morphology of block copolymer melts}, Macromolecules, 19, 2621-2632, 1986.

\bibitem{Hamley} I. Hamley, \emph{Developments in block copolymer science and technology}, Wiley, New York, 2004.

\bibitem{Bats_Fredrickson} F. S. Bats, G. H. Fredrickson, \emph{Block copolymers - designer soft materials}, Phys. Today 52 (2) (1999) 32.

\bibitem{Choksi} R. Choksi, \emph{Scaling laws in microphase separation of diblock copolymers}, J. Nonlinear Sci. 11 (2011) 223-236.

\bibitem{Bahiana_Oono} M. Bahiana, Y.Oono, \emph{Cell dynamical system approach to block copolymers}, Phys. Rev. A 41 (1990) 6763.

\bibitem{Hasegawa_Tannaka_Yamasaki_Hashimoto} H. Hasegawa, H. Tannaka, K. Yamasaki, T. Hashimoto, \emph{Bicontinuous microdomain morphology of block copolymers. 1. tetrapod-network structure of polysrtrene-polyisoprene diblock polymers}, Macromolecules 20 (1987) 1651-1662.

\bibitem{Zheng_Wang} W. Zheng, Z.-G. Wang, \emph{Morphology of ABC triblock copolymers}, Macromolecules 28 (21) (1995) 7215-7223.

\bibitem{Tang_Qiu_Zhang_Yang} P. Tang, F. Qiu, H. Zhang, Y. Yang, \emph{Morphology and phase diagram of complex block copolymers: abc linear triblock copolymers}, Phys. Rev. E 69 (2004) 031803.

\bibitem{Li_Jiang_Chen} S. Li, Y. Jiang, J. Chen, \emph{Morphologies and phase diagrams of abc star triblock copolymers confined in a spherical cavity}, Macromolecules 9 (2013) 4843-4854.

\bibitem{Zhao_Xu2019} Y. Zhao, X. Xu, \emph{Energy Stable Semi-implicit Schemes for
Allen–Cahn–Ohta–KawasakiModel in Binary System}, J. Sci. Comput., 2019.

\bibitem{Zhao_Xu2020} Y. Zhao, X. Xu, \emph{Maximum Principle Preserving Schemes for Binary Systems with Long-Range Interactions}, J. Sci. Comput., 2020.

\bibitem{Du_Yang2016} Q. Du, J. Yang, \emph{Asymptotically compatible Fourier spectral approximations of nonlocal Allen-Cahn equations}, SIAM J. Numer. Anal., 54 (3), 1899–1919, 2016.

\bibitem{Du_Yang2017} Q. Du, J. Yang, \emph{Fast and accurate implementation of Fourier spectral approximations of nonlocal diffusion operators and its applications}, J. Comput. Phys., 332 (2017) 118–134.

\bibitem{Du_nonlocalbook} Q. Du, \emph{Nonlocal modeling, analysis and computation}, CBMS-NSF Regional Conference Series in Applied Mathematics, 94, 2020.

\bibitem{ShenYang_DCDS2010} J. Shen and X. Yang, \emph{Numerical approximations of Allen-Cahn and Cahn-Hillard equations}, Discrete Contin. Dyn. Syst. A, 28, 1669, 2010.

\bibitem{WiseWangLowengrub_SINA2009} S. Wise, C. Wang and J. Lowengrub, \emph{An energy stable and convergent finite difference scheme for the phase field crystal equation}, SIAM J. Numer. Anal. 47, 2269-2288, 2009.

\bibitem{HuWiseWangLowengrub_JCP2009} Z. Hu, S. Wise, C. Wang and J. Lowengrub, \emph{Stable and efficient finite-difference nonlinear-multigrid scheme for the phase field crystal equation}, J. Comput. Phys. 228, 5323-5339, 2009.

\bibitem{WangWise_SINA2011} C. Wang and S. Wise, \emph{An energy stable and convergent finite-difference scheme for the modified phase field crystal equation}, SIAM J. Numer. Anal. 49, 945-969, 2011.

\bibitem{ShenWangWangWise_SINA2012} J. Shen, C. Wang, X. Wang, S. Wise, \emph{Second-order convex splitting schemes for gradient flows with Ehrlich–Schwoebel type energy: application to think film epitaxy}, SIAM J. Numer. Anal. 50, 105–125 (2012).

\bibitem{ChenCondeWangWangWise_JSC2012} W. Chen, S. Conde, C. Wang, X. Wang and S. Wise, \emph{A linear energy stable scheme for a thin film model without slope selection}, J. Sci. Comput. 52, 546-562, 2012.

\bibitem{Eyre98} D. Eyre, \emph{Unconditionally gradient stable time marching the Cahn–Hillard equation}, Computational and Mathematical Models of Microstructural Evolution (San Francisco, CA, 1998), Materials Research Society Symposia Proceedings, vol. 529, p. 39 (1998).

\bibitem{BenesovaMelcherSuli_SJNA2014} B. Benesova, C. Melcher, E. Suli, \emph{An implicit midpoint spectral approximation of nonlocal Cahn-Hilliard equations}, SIAM J. Numer. Anal. 52 (2014) 1466.

\bibitem{XuTang_SJNA2006} C. Xu, T. Tang, \emph{Stability analysis of large time-stepping methods for epitaxial growth models}, SIAM J. Numer. Anal. 44 (2006) 1759.

\bibitem{Yang_JCP2016} X. Yang, \emph{Linear and unconditionally energy stable numerical schemes for the phase field model of homopolymer blends}, J. Comput. Phys. 302 (2016) 509.

\bibitem{ChengYangSHen_JCP2017} W. Cheng, X. Yang, J. Shen, \emph{Efficient and accurate numerical schemes for a hydro-dynamically coupled phase field diblock copolymer model}, J. Comput. Phys. 341, 44 (2017)

\bibitem{CaffarelliMuler_ARMA1995} L. Caffarelli, N. E. Muler, \emph{A L1 bound for solutions of the cahn-hilliard equation}, Arch. Rational. Mech. Anal. 133 (1995) 129-144.

\bibitem{DuJuLiQiao_JCP2018} Q. Du, L. Ju, X. Li, Z. Qiao, \emph{Stabilized linear semi-implicit schemes for the nonlocal Cahn–Hilliard equation}, J. Comput. Phys. 363 (2018) 39-54.

\bibitem{DuJuLiQiao_SINA2019} Q. Du, L. Ju, X. Li, Z. Qiao, \emph{Maximum Principle Preserving Exponential Time Differencing Schemes for the Nonlocal Allen-Cahn Equation}, SIAM J. Numer. Anal. 57 (2019) 875-898.

\bibitem{Zhao_Luo1d} Y. Zhao, W. Luo, \emph{Nonlocal Effects on a 1D Generalized Ohta-Kawasaki Model}, To be appeared in Physica D.

\bibitem{Zhao_Luo2d} Y. Zhao, W. Luo, \emph{Nonlocal effects on a 2D Generalized Ohta-Kawasaki model}, preprint.

\bibitem{Zhao_Luodisk} Y. Zhao, W. Luo, \emph{ACOK/ACNO Systems on the Unit Disk Domain}, preprint.

\bibitem{Zhao_Choi} Y. Zhao, HJ. Choi, \emph{Second-order Stabilized Semi-implicit Energy Stable Schemes for Bubble Assemblies in Binary and Ternary Systems}, Discrete Contin. Dyn. Syst. Ser. B, 2021.

\bibitem{Shen_Tang_Wang} J. Shen, T. Tang and L. Wang, \emph{Spectral Methods:Algorithms, Analysis and Applicatiions}, Springer Series in Computational Mathematics, 41. Springer, Heidelberg, 2011. 

\bibitem{TianDu_SINA} X. Tian and Q. Du, \emph{Asymptotically compatible schemes and applications to robust discretization of nonlocal models}, SIAM J. Numer. Anal., 52(2014), pp. 1641–1665. 

\bibitem{XuRussellOckoChecco_SoftMatter2011} J. Xu, T. Russell, B. Ocko, A. Checco, \emph{Block copolymer self-assembly in chemically patterned squares}, Soft matter 7 (2011) 3915.

\bibitem{Nishiura_Ohnishi1995} Y. Nishiura and I. Ohnishi, \emph{Some mathematical aspects of the microphase separation in diblock copolymers}, Physica D, 84(1-2):31–39, 1995.

\bibitem{Ren_Wei2000} X. Ren and J. Wei, \emph{On the multiplicity of solutions of two nonlocal variational problems}, SIAM J. Math. Anal., 31(4):909–924, 2000.

\bibitem{Choksi2012} R. Choksi, \emph{On global minimizers for a variational problem with long-range interactions}, Quart. Appl. Math. 70 (2012) 517–537.

\bibitem{Xu_Zhao2019} X. Xu and Y. Zhao, \emph{Energy Stable Semi-implicit Schemes for Allen–Cahn–Ohta–Kawasaki ModelinBinarySystem}, J. Sci. Comput., 80(2019), 1656-1680. 

\bibitem{Xu_Zhao2020} X. Xu and Y. Zhao, \emph{Maximum principle preserving schemes for binary systems with long-range interaction}, J. Sci. Comput., 84(2020), 34pp.

\bibitem{Choi_Zhao2021} H. Choi and Y. Zhao, \emph{Second-order stabilized semi-implicit energy stable schemes for bubble assemblies in binary and ternary systems}, Discrete Contin. Dyn. Syst. Ser. B, 2021.

\bibitem{DuGunzburgerLehoucqZhou_SIRE2012} Q. Du, M. Gunzburger, R. Lehoucq, K. Zhou, \emph{Analysis and approximation of nonlocal diffusion problems with volume constraints}, SIAM Rev. 54 (2012) 667–696.

\bibitem{DuGunzburgerLehoucqZhou_M3AS2013} Q. Du, M. Gunzburger, R. Lehoucq, K. Zhou, \emph{A nonlocal vector calculus, nonlocal volume-constrained problems, and nonlocal balance laws}, Math. Models Methods Appl. Sci., 23 (2013), pp. 493–540.

\bibitem{DuZhou_2013} Q. Du and K. Zhou, \emph{Mathematical analysis for the peridynamic nonlocal continuum theory}, Math. Model. Numer. Anal., 45 (2011), pp. 217–234.

\bibitem{BurchLehoucq_2011} N. Burch and R. B. Lehoucq, \emph{Classical, nonlocal, and fractional diffusion equations on
bounded domains}, Internat. J. Multiscale Comput. Engrg., 9 (2011), pp. 661–674.

\bibitem{AndreuMazonRossiToledo_2010} F. Andreu, J. M. Mazon, J. D. Rossi, and J. Toledo, \emph{Nonlocal Diffusion Problems}, Math. Surveys Monographs 165, AMS, Providence, RI, 2010.

\bibitem{OhtaNakazawa_Macromolecules1993} H. Nakazawa, T. Ohta, \emph{Microphase separation of ABC-type triblock copolymers}, Macromolecules 26 (1993) 5503–5511.

\bibitem{RenWei_JNS2003} X. Ren, J. Wei, \emph{Triblock copolymer theory: Ordered ABC lamellar phase}, J. Nonlinear Sci. 13 (2003)
175–208.

\bibitem{RenWei_ARMA2013} X. Ren, J. Wei, \emph{A double bubble in a ternary system with inhibitory long range interaction}, Arch. Ration. Mech. Anal. 208 (2013) 201–253.

\bibitem{RenWei_ARMA2015} X. Ren, J. Wei, \emph{A double bubble assembly as a new phase of a ternary inhibitory system}, Arch. Rat. Mech. Anal. 215 (3) (2015) 967–1034.

\bibitem{GennipPeletier_CVPDE2008} Y. Gennip, M. Peletier, \emph{Copolymer-homopolymer blends: global energy minimisation and global energy bounds}, Calc. Var. Partial Differ. Equ. 33 (2008) 75–111.

\bibitem{XuDu_JNS2022}Z. Xu, Q. Du, \emph{On the ternary ohta-kawasaki free energy and its one dimensional global minimizers}, J Nonlinear Sci 32, 61 (2022).

\bibitem{JJ_NC2020} J. Jena. et al., \emph{Elliptical bloch skyrmion chiral twins in an antiskyrmion system}, Nature Commu. 11 (2020) 1115.

\bibitem{Wang_RenCCM2019} X. Ren and C. Wang, \emph{Stationary disk assemblies in a ternary system with long range interaction}, Communications in Contemporary Mathematics, Vol 21, No. 06 (2019) 1850046.

\bibitem{Du_WangJMB2008} Q. Du and X. Wang, \emph{Modelling and Simulations of Multi-component Lipid Membranes and Open Membranes via Diffuse Interface Approaches}, J. Mathematical Biology, 56, no3, 347-371, 2008


\bibitem{DuTian_FCM2020} Q. Du,  X. Tian, \emph{Mathematics of Smoothed Particle Hydrodynamics: A Study via Nonlocal Stokes Equations}, Foundations of Computational Mathematics, (2020) 20:801–826

\end{thebibliography}

\end{document}